\documentclass[12pt]{amsart}
\usepackage{amsmath,amsfonts,euscript,amscd,amsthm,amssymb,upref,graphics}
\usepackage[all]{xy}
\usepackage{color}

\theoremstyle{definition}

\swapnumbers

\theoremstyle{plain}

\newtheorem{theorem}{Theorem}[section]
\newtheorem{slicetheorem}[theorem]{Slice Theorem}

\newtheorem{proposition}[theorem]{Proposition}
\newtheorem{lemma}[theorem]{Lemma}

\newtheorem{corollary}[theorem]{Corollary}

\newtheorem{nothing}[theorem]{}

\theoremstyle{definition}

\newtheorem{definition}[theorem]{Definition}

\newtheorem{nothing*}[theorem]{}
\newtheorem{example}[theorem]{Example}

\newtheorem{examples}[theorem]{Examples}

\newtheorem{notation}[theorem]{Notation}

\newtheorem{subnothing*}[sub]{}

\newtheorem{bigremark}[theorem]{Remark}
\newtheorem{bigremarks}[theorem]{Remarks}

\newtheorem{caution}[theorem]{Caution}

\theoremstyle{remark}

{\begin{enumerate}\setlength{\itemsep}{#1}}{\end{enumerate}}

\newcommand{\Spec}{	\operatorname{\text{\rm Spec}}}

\newcommand{\haut}{	\operatorname{\text{\rm ht}}}
\newcommand{\supp}{	\operatorname{\text{\rm supp}}}

\newcommand{\trdeg}{	\operatorname{\text{\rm trdeg}}}
\newcommand{\Frac}{	\operatorname{\text{\rm Frac}}}
\newcommand{\HFrac}{	\operatorname{\text{\rm HFrac}}}
\newcommand{\Char}{	\operatorname{\text{\rm char}}}

\newcommand{\lnd}{	\operatorname{\text{\rm LND}}}
\newcommand{\hlnd}{	\operatorname{\text{\rm HLND}}}
\newcommand{\Gr}{	\operatorname{\text{\rm Gr}}}

\newcommand{\Proj}{	\operatorname{\text{\rm Proj}}}
\newcommand{\Div}{	\operatorname{\text{\rm Div}}}
\renewcommand{\div}{	\operatorname{{\rm div}}}
\newcommand{\lcm}{	\operatorname{\text{\rm lcm}}}

\newcommand{\EXT}{	\operatorname{\text{\rm EXT}}}
\newcommand{\HT}{	\operatorname{\text{\rm HT}}}

\newcommand{\YY}{Y^{(1)}}

\newcommand{\PP}{\operatorname{{\rm\bf P}_\bk}}
\newcommand{\GND}{\operatorname{{\rm\bf GND}_\bk}}

\newcommand{\setspec}[2]{\big\{\,#1\, \mid \,#2\, \big\}}

\newcommand{\Integ}{\ensuremath{\mathbb{Z}}}
\newcommand{\Nat}{\ensuremath{\mathbb{N}}}
\newcommand{\Rat}{\ensuremath{\mathbb{Q}}}

\newcommand{\aff}{\ensuremath{\mathbb{A}}}
\newcommand{\proj}{\ensuremath{\mathbb{P}}}
\newcommand{\bk}{{\ensuremath{\rm \bf k}}}

\newcommand{\bbV}{\ensuremath{\mathbb{V}}}
\newcommand{\bbD}{\ensuremath{\mathbb{D}}}

\newcommand{\pgoth}{{\ensuremath{\mathfrak{p}}}}
\newcommand{\qgoth}{{\ensuremath{\mathfrak{q}}}}
\newcommand{\Pgoth}{{\ensuremath{\mathfrak{P}}}}

\newcommand{\Aeul}{\EuScript{A}}

\newcommand{\Leul}{\EuScript{L}}

\newcommand{\Oeul}{\EuScript{O}}

\newcommand{\Reul}{\EuScript{R}}

\newcommand{\Xeul}{\EuScript{X}}

\newcommand{\isom}{\cong}
\renewcommand{\epsilon}{\varepsilon}
\renewcommand{\phi}{\varphi}
\renewcommand{\emptyset}{\varnothing}

\newenvironment{enumerata}%
{\begin{enumerate}

}{\end{enumerate}}

\newcommand{\rien}[1]{}

\begin{document}
\renewcommand{\baselinestretch}{1.07}


\title[Locally nilpotent derivations of graded integral domains]{Locally nilpotent derivations of\\ graded integral domains and cylindricity}

\author{Michael Chitayat and Daniel Daigle}

\address{Department of Mathematics and Statistics\\
University of Ottawa\\
Ottawa, Canada\ \ K1N 6N5}

\email{mchit007@uottawa.ca}
\email{ddaigle@uottawa.ca}

\thanks{Research of both authors supported by grant 04539/RGPIN/2015 from NSERC Canada.}


{\renewcommand{\thefootnote}{}
}

\begin{abstract} 
Let $B$ be a commutative $\Integ$-graded domain of characteristic zero.
An element $f$ of $B$ is said to be {\it cylindrical\/} if it is nonzero, homogeneous of nonzero degree,
and such that $B_{(f)}$ is a polynomial ring in one variable over a subring.
We study the relation between the existence of a cylindrical element of $B$ and the existence of a nonzero locally nilpotent derivation of $B$.
Also, given $d\ge1$, we give sufficient conditions that guarantee that every derivation of $B^{(d)} = \bigoplus_{i \in \Integ} B_{di}$ can be extended to 
a derivation of $B$.
We generalize some results of Kishimoto, Prokhorov and Zaidenberg that relate
the cylindricity of a polarized projective variety $(Y,H)$ to the existence of a nontrivial $G_a$-action on the affine cone over $(Y,H)$.
\end{abstract}
\maketitle
  
\vfuzz=2pt


\section{Introduction}

Let $B$ be a ring (by which we mean a commutative, associative, unital ring).
A derivation $D : B \to B$ is said to be {\it locally nilpotent} if, for each $b \in B$, there exists $n>0$ such that $D^n(b)=0$.
We write $\lnd(B)$ for the set of locally nilpotent derivations $D : B \to B$.
The ring $B$ is said to be {\it rigid\/} if $\lnd(B) = \{0\}$.

It has been known for a long time that the rigidity of $B$ is related to the geometry of $\Spec B$ in the following way.
If $B$ is an affine domain over a field $\bk$ of characteristic zero,
then $B$ is non-rigid if and only if there exists $f \in B \setminus \{0\}$ such that the basic open set
$\bbD(f) = \setspec{\pgoth \in \Spec B}{ f \notin \pgoth }$
is isomorphic to $\aff^1_\bk \times Z$ for some affine variety~$Z$.

A graded ring $B$ is said to be {\it rigid\/} if $\lnd(B) = \{0\}$, i.e., if it is rigid as a non-graded ring.
If $B$ is $\Nat$-graded, it is natural to ask whether the rigidity of $B$ is related to the geometry of the scheme $\Proj B$,
and in particular to the existence of an open subset $\bbD_+(f) = \setspec{\pgoth \in \Proj B}{ f \notin \pgoth }$
of  $\Proj B$ isomorphic to a product  $\aff^1_\bk \times Z$.
To discuss this question, it is convenient to introduce the following definition:
\begin{quote}
An element $f$ of a $\Integ$-graded domain $B$ is {\it cylindrical\/} if it is nonzero, homogeneous of nonzero degree,
and is such that $B_{(f)}$ is a polynomial ring in one variable.
\end{quote}
By $B_{(f)}$, we mean the degree-$0$ subring of the graded ring $B_f = S^{-1}B$ where $S= \{1,f,f^2,\dots\}$;
saying that $B_{(f)}$ is a polynomial ring in one variable means that there exists a subring $R$ of $B_{(f)}$ such
that $B_{(f)}$ is a polynomial ring in one variable over $R$.

Now, let $B$ be an $\Nat$-graded affine domain over a field $\bk$ of characteristic zero.
Clearly, $B$ has a cylindrical element if and only if
there exists  a homogeneous $f \in B\setminus\{0\}$ of positive degree such that 
$\bbD_+(f) \isom \aff^1_\bk \times Z$ for some affine variety $Z$.
So, returning to the question of how the rigidity of $B$ is related to the geometry of $\Proj B$,
one natural question to ask is the following:
{\it how is the rigidity of $B$ related to the existence of a cylindrical element of $B$?}
Note that this question is formulated in purely algebraic language.

Theorem 0.6 and Corollary 3.2 of \cite{KPZ2013} answer this question in the context of the Dolgachev-Pinkham-Demazure (DPD) construction
(these two results of \cite{KPZ2013} are quoted in Section \ref{jhfbo293wei0dc02uew9}).
The DPD construction allows us to go back and forth between $\Nat$-graded normal affine $\bk$-domains $B$ and polarized projective $\bk$-varieties $(Y,H)$.
The results of \cite{KPZ2013} relate the rigidity of $B$ to the geometric properties of the associated polarized projective variety $(Y,H)$.
The fact that the results of \cite{KPZ2013} are written in the formalism of the DPD construction has benefits and disadvantages. 
The obvious benefit is that the body of knowledge concerning polarized projective varieties
can be brought to bear for studying the rigidity of graded rings.
The disadvantage is that those two results of \cite{KPZ2013} are difficult to use in an algebraic setting.

The aim of the present article is to give algebraic analogues and generalizations of the results of \cite{KPZ2013}.
In doing so, we fill a gap in the proof of Theorem 0.6 of \cite{KPZ2013} (see Remark \ref{aPOKJxcbcdgfjkawEv3mnCigwvgdEduh9w} for details).
We also give some results on the extension of derivations.

Recall that the {\it saturation index\/} $e(B)$ of a graded domain $B = \bigoplus_{i \in \Integ} B_i$ 
is defined by $e(B) = \gcd\setspec{ i \in \Integ }{ B_i \neq 0 }$.
We claim that Theorem 0.6 and Corollary 3.2 of \cite{KPZ2013}, once translated into algebraic language,
are equivalent to parts (1) and (2), respectively, of the following assertion:

\begin{nothing}  \label {iocjvn2uygew8dhc}
\it Let $B = \bigoplus_{i \in \Nat} B_i$ be an $\Nat$-graded normal domain that is finitely generated as an algebra over an algebraically closed
field of characteristic zero, and assume that the transcendence degree of $B$ over $B_0$ is at least $2$ and that $e(B)=1$.
\begin{enumerate}

\item The following hold:
\begin{enumerate}

\item If $B$ is non-rigid, then $B$ has a cylindrical element.

\item If $B$ has a cylindrical element then there exists $d\ge1$ such that $B^{(d)}$ is non-rigid,
where $B^{(d)} = \bigoplus_{ i \in \Nat } B_{id}$.

\end{enumerate}

\item If $e(B/\pgoth) = 1$ for all $\pgoth \in \Proj(B)$, then
$$
\text{$B$ is non-rigid if and only if $B$ has a cylindrical element.}
$$
\end{enumerate}
\end{nothing}

The claim that \ref{iocjvn2uygew8dhc} is equivalent to the two aforementioned results of \cite{KPZ2013} is not entirely obvious.
Section \ref{jhfbo293wei0dc02uew9} is devoted to proving this claim.

We say that $B$ is {\it saturated in codimension~$1$}
if $e(B/\pgoth)=e(B)$ for every homogeneous prime ideal $\pgoth$ of $B$ of height $1$.
The following generalizes \ref{iocjvn2uygew8dhc} (and hence Theorem 0.6 and Corollary 3.2 of \cite{KPZ2013}):

\begin{theorem} \label {edh83yf6r79hvujhxu6wrefji9e} 
Let $B = \bigoplus_{i \in \Integ}B_i$ be a $\Integ$-graded domain that satisfies one of:
\begin{itemize}

\item[(i)] the transcendence degree of $B$ over $B_0$ is at least $2$;

\item[(ii)] there exist $i,j$ such that $i < 0 < j$, $B_i \neq 0$ and $B_j \neq 0$.

\end{itemize}
Assume that $B$ is finitely generated as an algebra over a field of characteristic zero.
\begin{enumerate}

\item The following are equivalent:
\begin{enumerate}

\item There exists $d \ge 1$ such that $B^{(d)}$ is non-rigid,
where $B^{(d)} = \bigoplus_{ i \in \Integ } B_{id}$.

\item $B$ has a cylindrical element.

\end{enumerate}

\item If $B$ is saturated in codimension $1$ and $B$ is normal, then conditions {\rm (a)} and {\rm (b)} are equivalent to:
\begin{enumerate}

\item[(c)] $B$ is non-rigid.

\end{enumerate}
\end{enumerate}
\end{theorem}

The proof of Theorem \ref{edh83yf6r79hvujhxu6wrefji9e} can be found in paragraph \ref{ajsdfo93ew0sdpw0edp}.
Note that, in both (1) and (2), our assumptions on the grading and the base field are weaker than in \ref{iocjvn2uygew8dhc}.
In part (1), we are not assuming that $B$ is normal.
In part (2), the hypothesis ``$B$ is saturated in codimension $1$''
is considerably weaker than the assumption  ``$e(B/\pgoth) = 1$ for all $\pgoth \in \Proj(B)$'' of \ref{iocjvn2uygew8dhc}(2)
(to get an idea of how much generality is gained, see Example \ref{kjcv9w9efjwW0emc7qd}).
Note that  Theorem \ref{edh83yf6r79hvujhxu6wrefji9e} is a corollary of more general and more precise results that can be found in Section \ref{doicvb923w9hpdc09w}.

The following simple observation (proved in \ref{cjOlKo10eof9ccnqpkeuycazvds4wud}) also deserves to be mentioned:

\begin{proposition}  \label {kcnboiqwdcm203}
Let $B$ be a $\Integ$-graded domain containing a field $\bk$.
If $B$ has a cylindrical element then there exists a field $K$ such that $\bk \subseteq K \subseteq \Frac(B)$ and $\Frac(B) = K^{(2)}$.
\end{proposition}

Some of our proofs in Section \ref{doicvb923w9hpdc09w} were obtained from those of \cite{KPZ2013}
by a process that could be described as ``removing the geometry''.
Once the geometry is removed, one realizes that some hypotheses are no longer necessary.
Except for Section 5, which mixes algebra and geometry, our proofs are completely algebraic.

Section \ref{di864f2gdhj123n8erhd8} studies the question whether derivations  of $B^{(d)}$ can be extended to derivations of $B$,
where $B$ is a $\Integ$-graded ring.
The following is a special case of Theorem~\ref{kjchvtdes334eOF5A5WF6I4ut9ity0o6plkg}:

\begin{corollary}  \label {eocnb9nxrd4g5sywjaafgbb}
Let $B = \bigoplus_{i \in \Integ} B_i$ be a $\Integ$-graded noetherian normal domain containing $\Rat$.
If $B$ is saturated in codimension $1$ then, for every $d > 0$, the following hold:
\begin{enumerata}

\item Every derivation $\delta : B^{(d)} \to B^{(d)}$ extends uniquely to a derivation $D : B \to B$.
\item Every locally nilpotent derivation $\delta : B^{(d)} \to B^{(d)}$ extends uniquely to a locally nilpotent derivation $D : B \to B$.

\end{enumerata}
\end{corollary}

Note that part (2) of Theorem \ref{edh83yf6r79hvujhxu6wrefji9e} is an immediate consequence of part (1) and of Corollary \ref{eocnb9nxrd4g5sywjaafgbb}.
When $B$ is not assumed to be saturated in codimension $1$, Theorem~\ref{kjchvtdes334eOF5A5WF6I4ut9ity0o6plkg} gives information
about the values of $d$ for which all derivations of $B^{(d)}$ extend to derivations of $B$.
The material contained in \ref{jdnfo9w3upoefwe0pk}--\ref{kjcv9w9efjwW0emc7qd} facilitates the determination of those values of $d$ in applications.

\medskip

We thank M.\ Zaidenberg and Y.\ Prokhorov for answering our questions about the DPD construction.

\subsection*{Conventions}
If $B$ is an algebra over a ring $A$, the notation $B = A^{[n]}$ (where $n \in \Nat$) means that $B$ is isomorphic as an $A$-algebra to a polynomial
ring in $n$ variables over $A$.
If $L/K$ is a field extension, we write $L = K^{(n)}$ to indicate that $L$ is a purely transcendental extension of $K$ of transcendence degree $n$.
The word ``domain'' means ``integral domain'',
and ``affine $\bk$-domain'' (where $\bk$ is a field) means a domain which is a finitely generated $\bk$-algebra.
We write $\Frac B$ for the field of fractions of a domain $B$.
If $A \subseteq B$ are domains, $\trdeg(B:A)$ denotes the transcendence degree of $\Frac B$ over $\Frac A$.
We use ``$\setminus$'' for set difference, ``$\subset$'' for strict inclusion and  ``$\subseteq$'' for general inclusion.
We write $\Rat_{>0} = \setspec{ x \in \Rat }{ x>0 }$ and we follow the convention that $0 \in \Nat$.

\section{Preliminaries}

\begin{nothing*}
Let $G$ be an abelian group (with additive notation).
A {\it $G$-grading\/} of a ring $B$ is a family $\big( B_i \big)_{i \in G}$ of subgroups of $(B,+)$ satisfying $B = \bigoplus_{i \in G} B_i$
and $B_iB_j\subseteq B_{i+j}$ for all $i,j \in G$.
The phrase ``let  $B = \bigoplus_{i \in G} B_i$ be a $G$-graded ring'' means that we are considering the ring $B$ together with 
the $G$-grading $\big( B_i \big)_{i \in G}$.

Let  $B = \bigoplus_{i \in G} B_i$ be a $G$-graded ring.
\begin{enumerate}

\item An element of $B$ is {\it homogeneous\/} if it belongs to $\bigcup_{i \in G} B_i$.
If $x$ is a nonzero homogeneous element, the {\it degree\/} of $x$, $\deg(x)$, is the unique $i \in G$ such that $x \in B_i$.
The degree of a non-homogeneous element is not defined.

\item $B_0$ is a subring of $B$ and is called the {\it degree-$0$ subring\/} of $B$.

\item Given a homogeneous element $f$ of $B$, $B_{(f)}$ denotes the degree-$0$ subring of the $G$-graded ring $B_f = S^{-1}B$ where $S= \{1,f,f^2,\dots\}$.
Given a homogeneous prime ideal $\pgoth$ of $B$, $B_{(\pgoth)}$ denotes the degree-$0$ subring of the $G$-graded ring
$S^{-1}B$ where $S$ is the set of homogeneous elements of $B \setminus \pgoth$.

\item Given a subgroup $H$ of $G$,
we define the graded subring $B^{(H)}$ of $B$ by $B^{(H)} = \bigoplus_{i \in H} B_i$.
Given $d \in G$, let $\langle d \rangle$ denote the subgroup of $G$ generated by $d$;
then $B^{(\langle d \rangle)}$ is abbreviated to $B^{(d)}$, i.e., 
we define $B^{(d)} = \bigoplus_{i \in \langle d \rangle} B_i$.

\end{enumerate}
\end{nothing*}

The following is surely well known but we don't know a reference, so we provide a proof.

\begin{lemma}  \label {kcjviyc2e2w45w6d2o9kmgvczxwer}
Let $G$ be an abelian group, $B = \bigoplus_{i \in G} B_i$ a $G$-graded ring, and $R$ a subring of $B_0$.
If $B$ is finitely generated as an $R$-algebra then so is $B^{(H)}$ for every subgroup $H$ of $G$.
\end{lemma}

\begin{proof}
Choose a finite generating set $\{g_1, \dots, g_r\}$ of the $R$-algebra $B$ such that each $g_i$ is homogeneous;
let $d_i = \deg(g_i)$ for $i = 1, \dots, r$,
note that 
$$
M = \setspec{ (k_1,\dots,k_r) \in \Integ^r }{ \textstyle \sum_{i=1}^r k_i d_i \in H }
$$
is a subgroup of $\Integ^r$ and define $N = M \cap \Nat^r$. By Proposition 8.3 on page 59 of \cite{GrilletBook_2001},
$N$ is a finitely generated submonoid of $\Nat^r$.
Let $\{e_1, \dots, e_s\}$ be a finite generating set of $N$, where $e_i = (e_{i,1}, \dots, e_{i,r})$ for each $i$.
Let $h_i = g_1^{e_{i,1}} \cdots g_r^{e_{i,r}}$ for $i = 1, \dots, s$.
Then it is easily verified that $h_1, \dots, h_s$ generate the $R$-algebra $B^{(H)}$.
\end{proof}

\begin{lemma} \label {Jdfuw6rvf0fofpfihcETGje7}
Let $M$ be an additive submonoid of $\Integ$ and let $d \in \Integ \setminus \{0\}$.

If $\gcd\big( \gcd(M), d \big) = 1$, then there exists $m \in M$ such that $\gcd(m,d)=1$.
\end{lemma}

\begin{proof}
We may assume that $d \neq \pm1$, otherwise the claim is trivial.
Let $p_1, \dots, p_n$ ($n\ge1$) be the prime factors of $d$.
Let $I$ be the ideal of $\Integ$ generated by $M$.
Then $I \nsubseteq p_i\Integ$ for all $i \in \{1, \dots, n\}$.
By \cite[(1.B)]{Matsumura}, $I \nsubseteq \bigcup_{i=1}^n p_i\Integ$.
Pick $\xi \in I$ such that $\xi \notin \bigcup_{i=1}^n p_i\Integ$.
Then $\xi = x-y$ for some $x,y \in M$.
Let $q = p_1 \cdots p_n \ge2$ and define $m = \xi+qy$. Then $m \notin \bigcup_{i=1}^n p_i\Integ$ (so $\gcd(m,d)=1$)
and $m = x + (q-1)y \in M$.
\end{proof}

The following fact appeared in \cite[Section 1]{Nouaze-Gabriel_1967} and \cite[Prop.~2.1]{Wright:JacConj}.

\begin{slicetheorem} \label {SliceThm}
Let $B$ be a $\Rat$-algebra, $D \in \lnd(B)$ and $A = \ker(D)$.
If $s\in B$ satisfies $Ds=1$ then $B = A[s] = A^{[1]}$
and $D = \frac{d}{ds} : A[s] \to A[s]$.
\end{slicetheorem}

A subring $A$ of a domain $B$ is said to be {\it factorially closed in $B$} if for every $x,y \in B$, 
the condition $xy \in A \setminus\{0\}$ implies $x,y \in A$.
The following is well known (see for instance \cite[Sec.\ 1.4]{Freud:Book-new} or \cite[5.3, 7.5]{Dai:IntroLNDs2010}):

\begin{lemma} \label {dkjhf2983ed9eje9}
If $B$ is a domain of characteristic zero and $D \in \lnd(B) \setminus \{0\}$
then $\trdeg(B : \ker D) = 1$ and $\ker(D)$ is factorially closed (and hence algebraically closed) in $B$.
\end{lemma}

\begin{nothing*}
Let $G$ be an abelian group and $B = \bigoplus_{i\in G} B_i$ a graded ring.
A derivation $D : B \to B$ is {\it homogeneous\/}
if there exists $d \in G$ such that $D(B_i) \subseteq B_{i+d}$ holds for all $i \in G$;
if $D$ is homogeneous and nonzero then $d$ is unique, and is called the {\it degree\/} of $D$.
Let $\hlnd(B)$ be the set of homogeneous locally nilpotent derivations $D : B \to B$.
The graded ring $B$ is said to be {\it rigid\/} if $\lnd(B) = \{0\}$, i.e., if it is rigid as a non-graded ring.
Graded rings $B$ satisfying $\hlnd(B)=\{0\}$ and $\lnd(B) \neq\{0\}$ do exist, and  are not rigid;
see for instance Propositions 6.5 and 6.6 of \cite{DaiFreudMoser}.
However, the next result states that the conditions
$\hlnd(B)=\{0\}$ and $\lnd(B) = \{0\}$ are equivalent if we assume that $G$ is torsion-free and that $B$
is an affine domain over a field of characteristic zero.
\end{nothing*}

\begin{lemma}  \label {dkjcnbv293qeokwdncbow90}
Let $G$ be a torsion-free abelian group and $R = \bigoplus_{i \in G} R_i$ a $G$-graded domain
which is finitely generated as an algebra over a field of characteristic zero.
If $R$ is non-rigid then there exists a nonzero homogeneous locally nilpotent derivation $D : R \to R$.
\end{lemma}

\begin{proof}
(See \cite[1.7, 1.9]{Dai:TameWild} for details.)
Since $G$ is torsion-free, we can (and shall) endow it with a total order (cf.\ \cite[Prop.\ 1.1.7]{AndersonFeil}).
Using that order, define a degree function
$\deg : R \to G \cup \{-\infty\}$ by declaring that if $x = \sum_{i \in G} x_i \in R$, $x_i \in R_i$, then
$\deg(x) = \max\setspec{i \in G}{ x_i \neq 0 }$ if $x \neq 0$ and $\deg(x)=-\infty$ if $x=0$.
Then one defines the associated graded ring $\Gr(R) = \bigoplus_{i \in G} B_{\le i} / B_{<i}$,
where $B_{\le i} = \setspec{ x \in B }{ \deg(x)\le i }$ and $B_{< i} = \setspec{ x \in B }{ \deg(x) < i }$.
Let $D : R \to R$ be any nonzero  locally nilpotent derivation.
By Theorem 1.7(a) of \cite{Dai:TameWild}, $\deg D$ is defined; this means that there exists $d \in G$ satisfying
(i)~$\deg(Dx)-\deg(x) \le d$ for all $x \in R\setminus\{0\}$, and
(ii)~$\deg(Dx)-\deg(x) = d$ for some $x \in R\setminus\{0\}$.
As is well known, it then follows that an associated derivation $\Gr(D) : \Gr(R) \to \Gr(R)$ 
is defined, and that $\Gr(D)$ is nonzero, homogeneous and locally nilpotent. Since $\Gr(R) \isom R$, this proves the claim.
\end{proof}

We offer a slightly improved version of the Theorem of Vasconcelos \cite{VASC}:

\begin{lemma} \label {ldkjxhcvi82ewdj0wd}
Let $A \subseteq B$ be domains of characteristic zero such that $B$ is integral over $A$.
If $\delta : A \to A$ is a locally nilpotent derivation and $D : B \to B$ a derivation that extends $\delta$,
then $D$ is locally nilpotent.
\end{lemma}

\begin{proof}
The Theorem of Vasconcelos says that this is true under the additional assumption that $\Rat \subseteq A$.
Note that $\Integ \subseteq A$ and let $S = \Integ\setminus\{0\}$.
We have $\Rat \subseteq S^{-1}A \subseteq S^{-1}B$,  $S^{-1}B$ is integral over $S^{-1}A$,
$S^{-1}\delta :  S^{-1}A \to S^{-1}A$ is locally nilpotent, and
$S^{-1} D :  S^{-1}B \to S^{-1}B$ extends $S^{-1}\delta$; so, by 
the Theorem of Vasconcelos, $S^{-1} D$  is locally nilpotent.
Consequently, $D$ is locally nilpotent.
\end{proof}

\section{Extension of derivations}
\label {di864f2gdhj123n8erhd8}

Given a $\Integ$-graded domain $B = \bigoplus_{i \in \Integ} B_i$, one may ask which $d\ge1$ have the property that
every derivation\footnote{The derivations that we are considering are not assumed to be homogeneous.}
of $B^{(d)}$ extends to a derivation of $B$,
and what hypotheses on $B$ guarantee that all $d\ge1$ have that property.
The present section investigates these questions.

If $B = \bigoplus_{i \in \Integ} B_i$ is a $\Integ$-graded domain then
one defines the {\it saturation index\/} $e(B)$ of $B$ by $e(B) = \gcd\setspec{ i \in \Integ }{ B_i \neq 0 }$.
We say that $B$ is {\it saturated in codimension $1$} if $e(B/\pgoth)=e(B)$ for every homogeneous prime ideal $\pgoth$ of $B$ of height $1$.

\begin{notation} \label {kjxccrw578dywthd93ieugo9we}
Let $B = \bigoplus_{i \in \Integ} B_i$ be a $\Integ$-graded domain such that $e(B)=1$.
\begin{enumerate}

\item Let $\Pi(B)$ be the set of prime numbers $p$ satisfying 
$$
\text{$p \mid e(B/\pgoth)$ for some height $1$ homogeneous prime ideal $\pgoth$ of $B$.}
$$

\item Let $\Pi^*(B) = \setspec{ d \in \Nat \setminus \{0\} }{ \text{no element of $\Pi(B)$ divides $d$} }$.

\end{enumerate}

Note that $\Pi^*(B)$ is closed under multiplication and that the following are equivalent:
\begin{itemize}

\item $B$ is saturated in codimension $1$,
\item $\Pi(B) = \emptyset$,
\item $\Pi^*(B) = \Nat \setminus \{0\}$.
\end{itemize}
We will see in Lemma \ref{jdnfo9w3upoefwe0pk} that $\Pi(B)$ is a finite set under mild hypotheses on $B$.
\end{notation}

We state the main result of this section:

\begin{theorem}  \label {kjchvtdes334eOF5A5WF6I4ut9ity0o6plkg}
Let $B = \bigoplus_{i \in \Integ} B_i$ be a $\Integ$-graded noetherian normal domain of characteristic zero and such that $e(B)=1$.
For each $d \in \Pi^*(B)$ such that $d$ is a unit of $B$, the following hold.
\begin{enumerata}

\item Every derivation $\delta : B^{(d)} \to B^{(d)}$ extends uniquely to a derivation $D : B \to B$.
\item Every locally nilpotent derivation $\delta : B^{(d)} \to B^{(d)}$ extends uniquely to a locally nilpotent derivation $D : B \to B$.

\end{enumerata}
\end{theorem}

The proof of Theorem \ref{kjchvtdes334eOF5A5WF6I4ut9ity0o6plkg} consists of \ref{98w3gvddb9wudjcudj}--\ref{kjcUvi7wdcjnK9dcveoeih9e}.
In paragraphs \ref{98w3gvddb9wudjcudj}--\ref{nvo2983epdk02qkw0}, the rings are not necessarily graded.

\begin{notation} \label {98w3gvddb9wudjcudj}
We write $(A,B) \in \EXT$ as an abbreviation for:
\begin{quote}
$A$ is a ring, $B$ is an $A$-algebra, and for every derivation $\delta : A \to A$ there exists a unique derivation $D: B \to B$
that makes the following diagram commute:
$$
\xymatrix@R=12pt{
B \ar[r]^D & B \\
A \ar[u] \ar[r]^{\delta} & A \, . \ar[u]
}
$$
\end{quote}
\end{notation}

\begin{examples}
If $L/K$ is a separable algebraic field extension, then $(K,L) \in \EXT$.
If $A$ is a ring and $S$ is a multiplicative set of $A$, then $(A,S^{-1}A) \in \EXT$.
If $A$ is a ring and $B=\{0\}$ is the null $A$-algebra, then $(A,B) \in \EXT$.
\end{examples}

\begin{lemma} \label {dh2763g5d4fc5a4dsErf56w7}
Let $B$ be a noetherian normal domain and $A$ a subring of $B$.
Suppose that $(\Frac A , \Frac B ) \in \EXT$ and that
there exists a family $( f_i )_{i \in I}$ of elements of $A \setminus \{0\}$ satisfying:
\begin{itemize}

\item $(A_{f_i},B_{f_i}) \in \EXT$ for every $i \in I$;

\item no height $1$ prime ideal of $B$ contains all $f_i$.

\end{itemize}
Then $(A,B) \in \EXT$.
\end{lemma}

\begin{proof}
First note that if  $B$ is a field then $(A,B) \in \EXT$ easily follows from the assumption $(\Frac A, B) \in \EXT$.
So we may assume that $B$ is not a field. This implies that there exists a height $1$ prime ideal of $B$,
so $( f_i )_{i \in I}$  is not the empty family.

Let $\delta : A \to A$ be a derivation.
Then $\delta$ extends uniquely to a derivation $\delta' : \Frac A \to \Frac A$;
since  $(\Frac A , \Frac B ) \in \EXT$,
$\delta'$ extends uniquely to a derivation $\Delta : \Frac B \to \Frac B$.

Note that if $D_1$ and $D_2$ are derivations $B \to B$ that extend $\delta$, then
their extensions $D_1', D_2' : \Frac B \to \Frac B$ are extensions of $\delta'$ (because $D_1', D_2'$ extend $\delta$)
and hence satisfy $D_1' = \Delta = D_2'$.
So $\delta$ has at most one extension to a derivation $B \to B$.

For each $i \in I$, let $\delta_i : A_{f_i} \to A_{f_i}$ be the unique derivation that extends $\delta$;
since $(A_{f_i},B_{f_i}) \in \EXT$, $\delta_i$ extends uniquely to a derivation $\Delta_i : B_{f_i} \to B_{f_i}$.
Since the extension $\Delta_i' : \Frac B \to \Frac B$ of $\Delta_i$ is an extension of $\delta'$, we must have $\Delta_i'=\Delta$.
Consequently, $\Delta(B) \subseteq B_{f_i}$. This shows that $ \Delta(B) \subseteq \bigcap_{i \in I} B_{f_i} $.

Let $\Spec^1(B)$ denote the set of height $1$ prime ideals of $B$.
If $\pgoth \in \Spec^1(B)$ then there exists $i \in I$ such that $f_i \notin \pgoth$; then $\Delta(B) \subseteq B_{f_i} \subseteq B_\pgoth$,
showing that 
$$
\textstyle
\Delta(B) \subseteq \bigcap_{\pgoth \in \Spec^1(B)} B_\pgoth = B,
$$
the last equality because $B$ is a noetherian normal domain.
Since $\Delta(B) \subseteq B$, the restriction $D : B \to B$ of $\Delta$ is a derivation that extends $\delta$.
So $(A,B) \in \EXT$.
\end{proof}

It is well known (cf.\ \cite{MasudaMiyanishi:Lifting09}) that locally nilpotent derivations can be lifted through an \'etale ring extension.
The following shows that the same is true for arbitrary derivations:

\begin{lemma}  \label {nvo2983epdk02qkw0}
Let $A$ be a ring, $A[\mathbf{X}] = A[X_1,\dots,X_n] = A^{[n]}$, $f_1,\dots,f_n \in A[\mathbf{X}]$, $B = A[\mathbf{X}]/(f_1,\dots,f_n)$,
and $\pi : A[\mathbf{X}] \to B$ the canonical homomorphism of the quotient ring.
Let $P \in A[\mathbf{X}]$ be the determinant of the Jacobian matrix $\frac{\partial(f_1,\dots,f_n)}{\partial(X_1,\dots,X_n)}$.
If $\pi(P)$ is a unit of $B$, then $(A,B) \in \EXT$.
\end{lemma}

\begin{proof}
Let $\delta: A \to A$ be a derivation.
Let $J^*$ denote the adjoint matrix of $J = \left( \frac{\partial f_i}{\partial X_j} \right)$;
then $J^*$ has entries in $A[\mathbf{X}]$ and $JJ^* = PI_n$.
Let $D_0 : A[\mathbf{X}] \to A[\mathbf{X}]$ be the derivation which extends $\delta$ and satisfies $D_0(X_j)=0$ for $1 \le j \le n$.
Define $g_1,\dots,g_n \in A[\mathbf{X}]$ by
$$
\left(\begin{smallmatrix} g_1 \\ \vdots \\ g_n \end{smallmatrix}\right) =
- J^* \, \left(\begin{smallmatrix} D_0(f_1) \\ \vdots \\ D_0(f_n) \end{smallmatrix}\right).
$$
Let $D_1 : A[\mathbf{X}] \to A[\mathbf{X}]$ be the unique $A$-derivation such that $D_1(X_j) = g_j$ for $1 \le j \le n$.
Then the derivation $D = P D_0 + D_1 : A[\mathbf{X}] \to A[\mathbf{X}]$
satisfies
$$
D(a) = P \delta(a)\ \text{(for $a\in A$)}
\qquad \text{and} \qquad D(X_i) = g_i \ \text{(for $1 \le i \le n$)}
$$
and a straightforward calculation shows that $D(f_i) = 0$ for all $i=1,\dots,n$.
So we may consider the derivation $\bar D =  D \pmod{(f_1,\dots,f_n)} : B\to B$.
Recall that $u = \pi(P)$ is a unit of $B$.
We have $\bar D(a) = u \delta(a)$ for all $a\in A$, so the derivation $u^{-1}\bar D : B \to B$ is an extension of $\delta$.

To prove uniqueness of extensions, it suffices to show that the only $A$-derivation of $B$ is $0$.
Suppose that $D : B \to B$ is an $A$-derivation. Let $J^{(\pi)} \in M_n(B)$ be the matrix whose $(i,j)$-th entry is
$\pi( \frac{\partial f_i}{\partial X_j} ) \in B$.
Then $\det\big(J^{(\pi)}\big) = \pi(\det J)$ is a unit of $B$, so $J^{(\pi)}$ is an invertible matrix.
Let $x_i = \pi(X_i)$ for $1 \le i \le n$.
For each $i \in \{1, \dots, n\}$ we have $0 = D(0) = D( f_i(x_1,\dots,x_n) )
= \sum_{j=1}^n \frac{\partial f_i}{\partial X_j}(x_1,\dots,x_n) D(x_j)
= \sum_{j=1}^n \pi( \frac{\partial f_i}{\partial X_j} ) D(x_j)$,
so 
$$
J^{(\pi)}  \, \left(\begin{smallmatrix} D(x_1) \\ \vdots \\ D(x_n) \end{smallmatrix}\right) 
= \left(\begin{smallmatrix} 0 \\ \vdots \\ 0 \end{smallmatrix}\right) \, .
$$
Since $J^{(\pi)}$ is invertible, we have $D(x_j)=0$ for all $j$, so $D=0$.
\end{proof}


\begin{notation}
Let $B = \bigoplus_{i \in \Integ} B_i$ be a $\Integ$-graded domain.
For each $d \in \Nat \setminus \{0,1\}$, let $\Xeul_d$ be the set of nonzero homogeneous $x \in B$ satisfying $\gcd( \deg(x), d ) = 1$.
\end{notation}

\begin{notation}
Let $S$ be a subset of a ring $R$.
We define the notation ``$\HT(S)>1$'' to mean that no height $1$ prime ideal of $R$ contains $S$.
(Caution: if $\pgoth$ is the zero ideal of a field $K$ then $\pgoth$ is a prime ideal of height $0$, so $\haut(\pgoth)=0$,
but $\HT(\pgoth)>1$ is also true.)
\end{notation}

\begin{lemma}  \label {kncfo923dkX0kcec9e}
Let $B = \bigoplus_{i \in \Integ} B_i$ be a $\Integ$-graded noetherian normal domain of characteristic zero such that $e(B)=1$.
Let $d \in \Nat \setminus \{0,1\}$ and assume that $d$ is a unit of $B$ and that $\HT(\Xeul_d)>1$.

Every derivation $\delta : B^{(d)} \to B^{(d)}$ extends uniquely to a derivation $D : B \to B$,
and if $\delta$ is locally nilpotent then so is $D$.
\end{lemma}

\begin{proof}
No height $1$ prime ideal of $B$ contains $\setspec{ x^d }{ x \in \Xeul_d }$.
Since $\Frac(B) / \Frac(B^{(d)})$ is an algebraic extension of fields of characteristic zero,
we have $(\Frac(B^{(d)}), \Frac(B)) \in \EXT$.
Let us prove:
\begin{equation} \label {jchvtrq54es5eescq}
\text{$\big( ( B^{(d)} )_{x^d} , B_x \big) \in \EXT$ for all $x \in \Xeul_d$.}
\end{equation}
If \eqref{jchvtrq54es5eescq} is true then $(B^{(d)},B) \in \EXT$ by Lemma \ref{dh2763g5d4fc5a4dsErf56w7}. 

Let $U$ be the set of nonzero homogeneous elements of $B^{(d)}$.
Then $\Aeul = U^{-1}B^{(d)} = \bigoplus_{i \in \Integ} \Aeul_i$ is a $\Integ$-graded domain, $e(\Aeul)=d$, $\Aeul_d \neq 0$ and $\Aeul_0$ is a field.
It is easy to see that, given any $y \in \Aeul_d \setminus \{0\}$, we have
$\Aeul = \Aeul_0[y^{\pm1}]$ where $y$ is transcendental over $\Aeul_0$.
In particular, $\Aeul$ is a normal domain. Let $K = \Frac\Aeul = \Frac( B^{(d)} )$.

Let $x \in \Xeul_d$.  We claim that
\begin{equation} \label {kjncp203ier9d}
\begin{minipage}[t]{.85\textwidth}
for every $r\in\Nat\setminus\{0,1\}$ such that $r \mid d$ and every $n \in \Integ\setminus\{0\}$,
no element $k$ of $K$ satisfies $nk^r = x^d$.
\end{minipage}
\end{equation}
Indeed, assume that $k \in K$ satisfies $nk^r = x^d$ and let $i = \deg(x)$.
Since $\Aeul_0$ is a field of characteristic zero, we have $\frac1n \in \Aeul_0$, so $k^r = \frac1n x^d \in \Aeul_{di}$,
so $k$ is integral over $\Aeul$.  Since $k \in \Frac \Aeul$ and $\Aeul$ is normal, we get $k \in \Aeul$.
Since $k \in \Aeul$ and $k^r=\frac1n x^d  \in \Aeul_{di}$ is homogeneous, it follows that $k$ is homogeneous and so $k \in \Aeul_{di/r}$.
We have $\gcd(i,d)=1$ because $x \in \Xeul_d$, so $r \nmid i$,
so $di/r \notin d\Integ$, so  $\Aeul_{di/r} = 0$; then $k=0$ and hence $x=0$, a contradiction. This proves \eqref{kjncp203ier9d}.
Now \eqref{kjncp203ier9d} together with Theorem 9.1 in Chapter VI of \cite{LangAlgebraText} implies that 
the polynomial $T^d - x^d \in K[T]$ is irreducible.  Thus,
\begin{equation}  \label {awe92983we09fd0}
\text{$T^d-x^d$ is the minimal polynomial of $x$ over $K$.}
\end{equation}

Let $A = (B^{(d)})_{x^d} \subseteq B_x$ and observe that $B_x$ is $\Integ$-graded and that $(B_x)^{(d)} = A$.
Let us argue that $A[x]=B_x$.  To see this, consider a nonzero homogeneous element $\xi$ of $B_x$.
Since $\deg(x)$ is a unit in $\Integ/d\Integ$, there exists $j$ such that $0 \le j < d$ and $\deg(\xi) \equiv j \deg(x) \pmod{d}$.
Then $\deg(\xi/x^j) \in d\Integ$, so $\xi/x^j \in (B_x)^{(d)} = A$, so $\xi \in A x^j \subseteq A[x]$. It follows that $A[x]=B_x$.

Consider the surjective $A$-homomorphism $\phi : A[T] \to B_x$ that sends $T$ to $x$.
We have $T^d-x^d \in \ker\phi$ and, by \eqref{awe92983we09fd0}, no polynomial in $A[T]$ of degree less than $d$ belongs to $\ker\phi$.
By the division algorithm, it follows that $\ker\phi = (T^d-x^d)$, so $A[T]/(T^d-x^d) \isom B_x$.
If $P \in A[T]$ denotes the determinant of the $1 \times 1$ Jacobian matrix $\frac{\partial(T^d-x^d)}{\partial T}$,
then $\phi(P) = dx^{d-1}$ is a unit in $B_x$ (because $d$ is a unit of $B$ by assumption).
So, by Lemma \ref{nvo2983epdk02qkw0}, we have $(A,B_x) \in \EXT$. This proves \eqref{jchvtrq54es5eescq};
as noted at the beginning of the proof, it follows that $(B^{(d)},B) \in \EXT$.

If $\delta \in \lnd(A)$ then, since $(B^{(d)},B) \in \EXT$,
there exists a derivation $D:B \to B$ that extends $\delta$.
Since $B$ is integral over $B^{(d)}$, Lemma \ref{ldkjxhcvi82ewdj0wd} implies that $D \in \lnd(B)$.
\end{proof}

\begin{lemma}  \label {dp98fho2039woej102i}
Let $B = \bigoplus_{i \in \Integ} B_i$ be a $\Integ$-graded domain such that $e(B)=1$.
For each $d \in \Nat \setminus \{0,1\}$ we have $\Xeul_d \neq \emptyset$ and the following are equivalent:
\begin{enumerata}

\item $\HT( \Xeul_d ) > 1$;

\item $\gcd( e( B/\pgoth ), d ) = 1$ for every height $1$ homogeneous prime ideal $\pgoth$ of $B$;

\item $d \in \Pi^*(B)$.

\end{enumerata}
\end{lemma}

\begin{proof}
Let $d \in \Nat \setminus \{0,1\}$.
Let $M=\setspec{ i \in \Integ }{ B_i \neq 0 }$.  Since $e(B)=1$, we have $\gcd(M)=1$, so $\gcd\big( \gcd(M),d\big)=1$.
Lemma \ref{Jdfuw6rvf0fofpfihcETGje7} then implies that there exists $m \in M$ such that  $\gcd(m,d)=1$.
Then $B_m \neq 0$ and $\emptyset \neq B_m \setminus \{0\} \subseteq \Xeul_d$, so  $\Xeul_d \neq \emptyset$. 

It is clear that (b) $\Leftrightarrow$ (c).  Let us prove (a) $\Leftrightarrow$ (b).

Assume that (a) holds and let $\pgoth$ be a height $1$ homogeneous prime ideal of $B$.
By (a), we have $\Xeul_d \nsubseteq \pgoth$; pick $x \in \Xeul_d \setminus \pgoth$ and let $m=\deg(x)$.
Then $\gcd(m,d)=1$ and $B_m \nsubseteq \pgoth$; the fact that  $B_m \nsubseteq \pgoth$ implies that $(B/\pgoth)_m \neq 0$, so $e( B/ \pgoth ) \mid m$,
so $\gcd( e( B/\pgoth ), d ) = 1$, showing that (b) holds.  So (a) implies (b).

We prove the converse by contradiction: assume that (b) is true and (a) is false.
Since (a) is false, there exists  a height $1$ prime ideal $\pgoth$ of $B$ such that $\Xeul_d \subseteq \pgoth$.
The ideal of $B$ generated by all homogeneous elements of $\pgoth$ is a homogeneous prime ideal $\qgoth$ of $B$ that contains $\Xeul_d$;
since $0 \notin \Xeul_d\neq\emptyset$, we have $\qgoth \neq 0$, so $0 \neq \qgoth \subseteq \pgoth$, so $\pgoth = \qgoth$, so $\pgoth$ is homogeneous.
By (b), it follows that $\gcd\big( e( B/\pgoth ), d \big) = 1$.
Consider the submonoid $M = \setspec{ i \in \Integ }{ (B/\pgoth)_i \neq 0 }$ of $\Integ$.
Since $\gcd\big( \gcd(M), d \big) = 1$, 
Lemma \ref{Jdfuw6rvf0fofpfihcETGje7} implies that there exists $m \in M$ such that $\gcd(m,d)=1$
(so $B_m \setminus \{0\} \subseteq \Xeul_d$) and $(B/\pgoth)_m \neq 0$,
so $B_m \nsubseteq \pgoth$, so $\Xeul_d \nsubseteq \pgoth$, a contradiction. So (b) implies (a).
\end{proof}

\begin{nothing*} \label {kjcUvi7wdcjnK9dcveoeih9e}
{\bf Proof of Theorem \ref{kjchvtdes334eOF5A5WF6I4ut9ity0o6plkg}.}
Let $d \in \Pi^*(B)$ and assume that $d$ is a unit of $B$.
If $d=1$ then clearly assertions (a) and (b) of Theorem \ref{kjchvtdes334eOF5A5WF6I4ut9ity0o6plkg} are true.
Assume that $d\ge2$. Then Lemma \ref{dp98fho2039woej102i} implies that $\HT( \Xeul_d ) > 1$,
so Lemma \ref{kncfo923dkX0kcec9e} implies that (a) and (b) hold.  \hfill\qedsymbol
\end{nothing*}

We now give some complementary information about $\Pi(B)$.

\begin{lemma} \label {jdnfo9w3upoefwe0pk}
Let $B = \bigoplus_{i \in \Integ} B_i$ be a $\Integ$-graded domain such that $e(B)=1$.
\begin{enumerata}

\item If $S$ is a subset of $\Integ\setminus\{0\}$ satisfying $\HT\big( \bigcup_{ i \in S } B_i \big) > 1$,
then each element of $\Pi(B)$ is a prime factor of some element of $S$.

\item If $\HT\big( \bigcup_{ i \in \Integ \setminus \{0\} } B_i \big) > 1$ and $B$ is noetherian, then $\Pi(B)$ is a finite set.

\item If the condition $\HT\big( \bigcup_{ i \in \Integ \setminus \{0\} } B_i \big) > 1$ is false, then $\Pi(B)$ is the set
of all prime numbers and $\Pi^*(B) = \{1\}$.

\end{enumerata}
\end{lemma}

\begin{proof}
(a) Suppose that $S \subseteq \Integ\setminus\{0\}$ satisfies $\HT\big( \bigcup_{ i \in S } B_i \big) > 1$, and let $p \in \Pi(B)$.
Then there exists a homogeneous prime ideal $\pgoth$ of $B$ of height $1$ such that $p \mid e(B/\pgoth)$.
Since $\HT\big( \bigcup_{ i \in S } B_i \big) > 1$, we have $\bigcup_{ i \in S } B_i \nsubseteq \pgoth$,
so there exists $i \in S$ such that $B_i \nsubseteq \pgoth$. Then $(B/\pgoth)_i \neq 0$, so $e(B/\pgoth) \mid i$, so $p \mid i$, which proves (a).

(b) Let $I$ be the ideal of $B$ generated by $\bigcup_{ i \in \Integ \setminus \{0\} } B_i$.
Since $B$ is noetherian, there exists a finite subset $F = \{f_1, \dots, f_r\}$ of $\bigcup_{ i \in \Integ \setminus \{0\} } (B_i\setminus\{0\})$
that generates $I$.
Let $d_i = \deg f_i$ ($1 \le i \le r$) and note that $S = \{ d_1, \dots, d_r \}$ is a subset of $\Integ \setminus \{0\}$
with the property that $\HT\big( \bigcup_{ i \in S } B_i \big) > 1$.
By part (a), each element of $\Pi(B)$ is a prime factor of some element of $S$. So $\Pi(B)$ is finite.

(c) If the condition $\HT\big( \bigcup_{ i \in \Integ \setminus \{0\} } B_i \big) > 1$ is false then some 
height $1$ prime ideal $\pgoth$ of $B$  satisfies $\bigcup_{ i \in \Integ \setminus \{0\} } B_i \subseteq \pgoth$.
Since $e(B)=1$, we have $B \neq B_0$,
so $\pgoth$ is a minimal prime over-ideal of $\bigcup_{ i \in \Integ \setminus \{0\} } B_i \nsubseteq \{0\}$,
so $\pgoth$ is homogeneous (the ideal generated by the homogeneous elements of $\pgoth$ is a nonzero prime ideal included in $\pgoth$,
so equal to $\pgoth$).
Then $\setspec{ i \in \Integ }{ (B/\pgoth)_i \neq 0 } = \{0\}$, so $e(B/\pgoth)=0$.
Since every prime number divides $0$, $\Pi(B)$ is the set of all prime numbers. Consequently, $\Pi^*(B) = \{1\}$.
\end{proof}

\begin{lemma} \label {jnf283yehd9jw937ry}
Let $B = \bigoplus_{i \in \Integ} B_i$ be a $\Integ$-graded noetherian domain such that $e(B)=1$,
and suppose that $B = B_0[x_1, \dots, x_n]$ where $n\ge2$ and, for each $i\in\{1,\dots,n\}$,
$x_i \neq 0$ is homogeneous of degree $d_i \in \Integ \setminus \{0\}$.
\begin{itemize}

\item Let $e_i = \gcd(d_1, \dots, \hat d_i , \dots, d_n)$ for each $i\in\{1,\dots,n\}$.
\item Let $U = \setspec{ i }{ \text{$1\le i \le n$ and $x_i$ is a unit of $B$} }$ and $U^c = \{1,\dots,n\} \setminus U$.
\item Let $E$ be the set of prime factors of $\prod_{i \in U^c} e_i$ ($E=\emptyset$ if $U^c=\emptyset$).
\end{itemize}
\begin{enumerata}
\item $E \subseteq \Pi(B)$
\item If $U \neq \emptyset$ then $\Pi(B) \subseteq$ set of  prime factors of $\gcd\setspec{d_i}{i \in U}$.
\item If $\HT\big( \bigcup_{ i \in \Integ \setminus \{0\} } B_i \big) > 1$ then $\Pi(B) \subseteq$ set of  prime factors of $d_1 \cdots d_n$.
\item If $\HT(\{x_i,x_j\})>1$ for every choice of distinct $i,j \in \{1,\dots,n\}$, then $\Pi(B)=E$.
\end{enumerata}
\end{lemma}

\begin{proof}
(a) Let $p \in E$. Then there exists $i \in U^c$ such that $p \mid e_i$.
Note that $e_i>1$.
For each $j \in \Integ \setminus e_i\Integ$, every monomial $a x_1^{m_1} \cdots x_n^{m_n}$ ($a \in B_0\setminus\{0\}$, $m_k \in \Nat$) belonging to $B_j$
must satisfy $m_i>0$. So $B_j \subseteq x_iB$ for each $j \in \Integ \setminus e_i\Integ$.
Since $B$ is noetherian and $x_i$ is neither $0$ nor a unit of $B$, there exists a height $1$ homogeneous prime ideal $\pgoth$ of $B$ satisfying $x_i \in \pgoth$.
Then $B_j \subseteq \pgoth$ for  each $j \in \Integ \setminus e_i\Integ$, i.e.,
$\setspec{ j \in \Integ }{ (B/\pgoth)_j \neq 0 } \subseteq e_i\Integ$.
Then $e_i \mid e(B/\pgoth)$, so all prime factors of $e_i$ belong to $\Pi(B)$, so $p \in \Pi(B)$.

(b) Let $i \in U$. Then no prime ideal of $B$ contains $B_{d_i}$, so in particular $\HT( B_{d_i} ) > 1$;
then Lemma \ref{jdnfo9w3upoefwe0pk}(a) applied to $S = \{ d_i \}$ shows that every element of $\Pi(B)$ divides $d_i$.

Assertion (c) follows by applying  Lemma \ref{jdnfo9w3upoefwe0pk}(a) to $S = \{ d_1, \dots, d_n \}$.

(d) By (a), it suffices to show that  $\Pi(B) \subseteq E$.
Let $p \in \Pi(B)$. For every choice of distinct $i,j \in \{1,\dots,n\}$, $p$ must divide one of $d_i, d_j$
(by Lemma \ref{jdnfo9w3upoefwe0pk}(a) applied to $S = \{ d_i, d_j \}$).
It follows that at most one $i \in \{1, \dots, n\}$ is such that $p \nmid d_i$.
So we can choose $i \in \{1,\dots,n\}$ such that $p \mid d_j$ for all $j \in \{1,\dots,n\} \setminus \{i\}$.
Then $p \mid e_i$.  If $i \in U$ then $p \mid d_i$ by (b), so $p$ divides $\gcd(d_i,e_i) = 1$, a contradiction.
This shows that $p \mid e_i$ for some $i \in U^c$, so $p \in E$.
\end{proof}

\begin{proposition} \label {8C733Br5er1A723rhd09}
Let $\bk$ be a field of characteristic zero and $B = \bigoplus_{i \in \Integ} B_i$ a $\Integ$-graded normal affine $\bk$-domain.
Suppose that $B = B_0[x_1, \dots, x_n]$ where $n\ge2$ and, for each $i\in\{1,\dots,n\}$,
$x_i$ is a homogeneous prime element of $B$ of degree $d_i \in \Integ \setminus \{0\}$.
Assume that 
\begin{itemize}

\item [(i)] $\gcd(d_1, \dots, \hat d_i , \dots, d_n)=1$ for every $i\in\{1,\dots,n\}$;

\item [(ii)] for every choice of distinct $i,j \in \{1,\dots,n\}$,  $x_i, x_j$ are not associates.

\end{itemize}
Then the following are equivalent:
\begin{enumerata}

\item $B$ is non-rigid;
\item $B^{(d)}$ is non-rigid for all $d\ge1$;
\item $B^{(d)}$ is non-rigid for some $d\ge1$.

\end{enumerata}
\end{proposition}

\begin{proof}
If $i,j \in \{1,\dots,n\}$ are distinct then, by assumption (ii),
$(x_i)$ and $(x_j)$ are distinct prime ideals of $B$ of height $1$, so $\HT(\{x_i,x_j\})>1$.
By Lemma \ref{jnf283yehd9jw937ry} together with assumption (i), it follows that $B$ is saturated in codimension $1$.
So (c) implies (a) by Theorem \ref{kjchvtdes334eOF5A5WF6I4ut9ity0o6plkg}.
It is obvious that (b) implies (c),
and the fact that (a) implies (b) can be seen as follows. Suppose that $B$ is non-rigid and consider $d\ge1$.
By Lemma \ref{dkjcnbv293qeokwdncbow90}, there exists a nonzero homogeneous locally nilpotent derivation $D : B \to B$.
Let $A = \ker(D)$.
We have $e(A)=1$ by Corollary 4.2 of \cite{DaiFreudMoser} together with assumptions (i) and (ii).
So we can pick a homogeneous $a \in A \setminus \{0\}$ such that the degree of the derivation $aD:B \to B$ is a multiple of $d$.
Then $aD$ maps $B^{(d)}$ into itself and, consequently, the restriction of $aD$ to $B^{(d)}$ is a nonzero element of $\lnd( B^{(d)} )$.
Thus, $B^{(d)}$ is non-rigid. This proves the claim.
\end{proof}

In view of the next result, it may be useful to recall that
if $R$ is an affine domain over a field then $R$ is a universally catenary noetherian domain.

\begin{corollary} \label {kjcnvo9wed0qmdnb7q}
Let $R$ be a universally catenary noetherian domain and $n\ge2$.
Define a $\Integ$-grading on $R[\mathbf{X}] = R[X_1, \dots, X_n] = R^{[n]}$
by choosing $d_1,\dots,d_n \in \Integ\setminus \{0\}$ such that $\gcd(d_1,\dots,d_n)=1$
and declaring that $R \subseteq R[\mathbf{X}]_0$ and that, for each $i\in\{1,\dots,n\}$, $X_i$ is homogeneous of degree $d_i$.
Suppose that $F$ is a homogeneous prime element of $R[\mathbf{X}]$ such that 
\begin{equation} \label {odkcvbo293wdj0w}
\text{$F \notin (X_i,X_j)$ for every choice of $i,j \in \{1, \dots, n\}$,}
\end{equation}
where $(X_i,X_j)$ is the ideal of $R[\mathbf{X}]$ generated by $X_i,X_j$.
Consider the $\Integ$-graded domain $B = R[\mathbf{X}]/(F)$.
Then $\Pi(B) = E$, where $E$ is defined in Lemma \ref{jnf283yehd9jw937ry}.
\end{corollary}

\begin{proof}
Let $\pi : R[\mathbf{X}] \to B$ be the canonical homomorphism of the quotient ring and define $x_i =\pi(X_i)$ for all $i \in \{1,\dots,n\}$.
Thus, $B = \bar R[x_1,\dots,x_n]$ where $\bar R = \pi(R)$.
Since $\bar R \subseteq B_0$, we have $B = B_0[x_1,\dots,x_n]$.
Note that $x_i\neq0$ for all $i$ (because of \eqref{odkcvbo293wdj0w}), so $e(B)=1$.
In view of Lemma \ref{jnf283yehd9jw937ry}, it suffices to check that $\HT(\{x_i,x_j\})>1$ for every choice of distinct $i,j \in \{1,\dots,n\}$.
Let $i,j$ be distinct elements of $\{1,\dots,n\}$.
Proceeding by contradiction, suppose that the condition $\HT(\{x_i,x_j\})>1$  is false;
then there exists a height $1$ prime ideal $\pgoth$ of $B$ such that $\{x_i,x_j\} \subseteq \pgoth$.
Let $\Pgoth = \pi^{-1}(\pgoth) \in \Spec R[\mathbf{X}]$.
Since $R$ is universally catenary, $R[\mathbf{X}]$ is catenary, so  all maximal prime chains $0=\qgoth_0 \subset \cdots \subset \qgoth_\ell = \Pgoth$
between $0$ and $\Pgoth$ in $R[\mathbf{X}]$ have the same length, equal to $\haut \Pgoth$.
Since $0 \subset (F) \subset \Pgoth$ is such a maximal prime chain, we have $\haut\Pgoth=2$.
Since $(X_i,X_j)$ is a prime ideal of $R[\mathbf{X}]$ of height $\ge2$ satisfying $(X_i,X_j) \subseteq \Pgoth$,
we have $\Pgoth = (X_i,X_j)$, so $F \in (X_i,X_j)$, contradicting \eqref{odkcvbo293wdj0w}.
This contradiction shows that $\HT(\{x_i,x_j\})>1$, which completes the proof.
\end{proof}

\begin{example}  \label {kjcv9w9efjwW0emc7qd}
Let $\bk$ be a field, $n\ge3$ and $a_1,\dots,a_n \in \Nat \setminus \{0\}$, and define $d_i = \lcm(a_1,\dots,a_n) / a_i$ ($1 \le i \le n$).
Note that $\gcd(d_1,\dots,d_n)=1$.
Define an $\Nat$-grading on $\bk[X_1,\dots,X_n]$ by declaring that, for each $i\in\{1,\dots,n\}$, $X_i$ is homogeneous of degree $d_i$.
Then $X_1^{a_1} + \cdots + X_n^{a_n}$ is a homogeneous prime element of $\bk[X_1,\dots,X_n]$,
so $B = \bk[X_1,\dots,X_n] / ( X_1^{a_1} + \cdots + X_n^{a_n} )$ is an $\Nat$-graded domain and $e(B)=1$.
Moreover, $B$ is normal, since $\Spec B$ is a hypersurface with at most one singular point.
One calls $B$ a {\it Pham-Brieskorn ring}.
Let $e_i = \gcd(d_1, \dots, \hat d_i , \dots, d_n)$ for each $i\in\{1,\dots,n\}$.
By Corollary \ref{kjcnvo9wed0qmdnb7q},
$$
\text{$\Pi(B)$ is the set of prime factors of $\textstyle \prod_{i=1}^n e_i$.}
$$
In particular, $B$ is saturated in codimension $1$ if and only if $e_1= \cdots = e_n = 1$.

Let us compare this with the condition ``$e(B/\pgoth) = 1$ for all $\pgoth \in \Proj(B)$'' of \ref{iocjvn2uygew8dhc}(2).
It is easy to see that, given distinct $i,j \in \{1,\dots,n\}$, there exists  $\pgoth_{i,j} \in \Proj(B)$ satisfying
$x_ix_j \notin \pgoth_{i,j}$ and $x_k \in \pgoth_{i,j}$ for all $k \in \{1,\dots,n\} \setminus \{i,j\}$;
moreover, for each $\pgoth \in \Proj(B)$ there exist distinct $i,j \in \{1,\dots,n\}$ such that  $x_ix_j \notin \pgoth$.
From these observations, it follows that  $B$ satisfies ``$e(B/\pgoth) = 1$ for all $\pgoth \in \Proj(B)$''
if and only if $d_1,\dots,d_n$ are pairwise relatively prime.
We may reformulate these conditions directly in terms of the tuple $(a_1,\dots,a_n)$ that defines
$B = \bk[X_1,\dots,X_n] / ( X_1^{a_1} + \cdots + X_n^{a_n} )$, as follows:
\begin{itemize}

\item $B$ is saturated in codimension $1$ if and only if 
$$
\text{$\lcm(a_1, \dots, \hat a_i, \dots, a_n) = \lcm(a_1,\dots,a_n)$ for every $i \in \{1,\dots,n\}$.}
$$

\item $B$ satisfies ``$e(B/\pgoth) = 1$ for all $\pgoth \in \Proj(B)$'' if and only if 
$$
\text{$\lcm(a_i,a_j) = \lcm(a_1,\dots,a_n)$ for every choice of distinct $i,j \in \{1,\dots,n\}$.}
$$

\end{itemize}
So, when applied to the class of Pham-Brieskorn rings, 
Theorem \ref{edh83yf6r79hvujhxu6wrefji9e}(2) is considerably more general than assertion \ref{iocjvn2uygew8dhc}(2).

Let us also point out that if $n\ge4$ and $B$ is saturated in codimension $1$ then, by Proposition \ref{8C733Br5er1A723rhd09},
the following are equivalent:
\begin{enumerata}

\item $B$ is non-rigid;
\item $B^{(d)}$ is non-rigid for all $d\ge1$;
\item $B^{(d)}$ is non-rigid for some $d\ge1$.

\end{enumerata}
\end{example}

\section{Cylindrical elements of a graded ring}
\label {doicvb923w9hpdc09w}

Given a ring $B$, the phrase ``$B$ is a polynomial ring in one variable'' means that there exists a subring $A$ of $B$
such that $B$ is a polynomial ring in one variable over $A$.
Observe that the zero ring is a polynomial ring in one variable.

\begin{definition} \label {skjdbco8w3uewndq0ws}
Let $G$ be an abelian group and $B$ a $G$-graded ring.
An element $f$ of $B$ is {\it cylindrical\/} if it is nonzero and homogeneous,
the element $\deg(f)$ of $G$ has infinite order, and $B_{(f)}$ is a polynomial ring in one variable.
\end{definition}

\begin{lemma}  \label {dpcinvbi293jhdbcnc38i}
Let $G$ be a torsion abelian group,
$S = \bigoplus_{i \in G} S_i$ a graded ring,
$R = \bigoplus_{i \in G} R_i$ a graded subring of $S$, and suppose that 
\begin{enumerate}

\item[(i)] $S = R[v] = R^{[1]}$ for some homogeneous element $v$ of $S$,

\item[(ii)] the set $H = \setspec{ i \in G }{ R_i \neq 0 }$ is a subgroup of $G$.

\end{enumerate}
Then one of the following conditions is satisfied:
\begin{enumerata}

\item $S_0 = (R_0)^{[1]}$,

\item there exists $e \in H \setminus \{0\}$ such that, if $k>0$ denotes the order of $e$ in $H$,
then, for every $r \in R_e \setminus\{0\}$, 
$(S_0)_{r^k}$ is a polynomial ring in one variable.

\end{enumerata}
\end{lemma}

\begin{proof}
Let $d = \deg(v) \in G$ and let $n \ge 1$ be the positive generator of the subgroup $\setspec{ i \in \Integ }{ id \in H }$ of $\Integ$.
We claim:
\begin{equation}  \label {d9cih8273wgdopkc}
S_0 \subseteq R[v^n] .
\end{equation}
Indeed, if $s \in S_0$ then $s = \sum_{i \in \Nat} a_i v^i$ for a unique family $(a_i)$ of elements of $R$.
Since $s$ and $v$ are homogeneous, and since $R[v]=R^{[1]}$, it follows that each $a_i$ is homogeneous and $\deg(a_iv^i) = 0$ for each $i$ such that $a_i\neq0$.
This last condition implies that $id \in H$  (and hence $n \mid i$) for each $i$ such that $a_i\neq0$, so $s \in R[v^n]$. This proves \eqref{d9cih8273wgdopkc}.

Note that if $nd=0$ then $\deg(v^n)=0$, so \eqref{d9cih8273wgdopkc} implies that $S_0 = R_0[v^n] = (R_0)^{[1]}$, so (a) is satisfied.

From now-on, assume that $nd \neq 0$. Define $e = nd \in H \setminus\{0\}$, let $k\ge1$ be the order of $e$ in $H$,
and note that $R_e \neq 0$ (since $R_i \neq 0$ for every $i \in H$).
Let $r \in R_e \setminus \{0\}$ and let us prove that $(S_0)_{r^k}$  is a polynomial ring in one variable.
Consider diagram (A):
\begin{equation} \label {odijbi28372wd9j}
\text{\rm(A)}\ \ \xymatrix{
S \ar@{=}[r] & R[v]  \\
S_0 \ar@{^{(}->}[u] \ar@{^{(}->}[r] & R[v^n] \ar@{^{(}->}[u]
} \qquad \quad
\text{\rm(B)}\ \ \xymatrix{
S_{r^k} \ar@{=}[r] & \Reul[v]  \\
(S_0)_{r^k} \ar@{^{(}->}[u] \ar@{^{(}->}[r] & \Reul[v^n] \ar@{^{(}->}[u] \ar@{=}[r] & \Reul[w] =\Reul^{[1]}
}
\end{equation}
and note that the four arrows in (A) are degree-preserving homomorphisms of graded rings.
Let $\Reul = R_{r^k} = \bigoplus_{i \in G} \Reul_i$.
Diagram (B) is obtained from diagram (A) by localization with respect to the multiplicative subset $\setspec{ (r^k)^i }{ i \in \Nat }$ of $R_0$;
so the arrows in (B) are injective degree-preserving homomorphisms of graded rings.
Note that $r$ is a unit of $\Reul$,
define $w = \frac{v^n}{r} \in \Reul[v]$ and observe that $\Reul[v^n] = \Reul[w] =\Reul^{[1]}$; see diagram (B).
Since $w$ is a homogeneous element of $\Reul[v]$ of degree $0$, we have
\begin{equation}  \label {icyt65dlfeu12iCdmc}
( \Reul[w] )_0 = \Reul_0[w] = (\Reul_0)^{[1]} .
\end{equation}
Note that $S_{r^k} = \big( \bigoplus_{i \in G} S_i \big)_{r^k} = \bigoplus_{i \in G} (S_i)_{r^k}$
and in particular $(S_{r^k})_0 = (S_0)_{r^k}$.
This means that $(S_0)_{r^k}$ is the degree-$0$ subring of the graded ring $S_{r^k} = \Reul[v]$.
So the image of $(S_0)_{r^k}$ in $\Reul[w]$ is precisely the degree-$0$ subring of $\Reul[w]$,
so $(S_0)_{r^k} = (\Reul_0)^{[1]}$ by \eqref{icyt65dlfeu12iCdmc}.
\end{proof}

\begin{proposition}  \label {ckjb20932gediqw0idAwk}
Let $B = \bigoplus_{i \in \Integ} B_i$ be a $\Integ$-graded $\Rat$-domain and $D : B \to B$ a nonzero homogeneous locally nilpotent derivation
such that $\ker D \nsubseteq B_0$.
Then some element of $D(B) \cap \ker(D)$ is a cylindrical element of $B$.
\end{proposition}

\begin{proof}
Let $A = \ker(D)$. Since $D(B) \cap A$  is a nonzero homogeneous ideal of $A$ and $A \nsubseteq B_0$,
we can pick a nonzero homogeneous element $h$ of $D(B) \cap A$ of nonzero degree. 
Let $m = \deg(h) \in \Integ\setminus\{0\}$, $G = \Integ/m\Integ$, and $\pi : \Integ \to G$ the canonical surjective homomorphism.
Define a $G$-grading $B = \bigoplus_{ j \in G } B_j'$ by declaring that $B_j' = \bigoplus_{i \in \pi^{-1}(j)} B_i$ for each $j \in G$.
Since $h-1$ is a $G$-homogeneous element of $B$, the ring $S = B/(h-1)$ is $G$-graded, say $S = \bigoplus_{j \in G} S_j$.
Let $\rho : B \to S$ be the canonical homomorphism of the quotient ring and note that $\rho(B_j') = S_j$ for all $j \in G$.
We note:
\begin{equation} \label {okcjbv8238iwf}
\text{the only $\Integ$-homogeneous element of $\ker\rho$ is $0$}
\end{equation}
(this follows from the fact that if $x,y \in B \setminus\{0\}$ are such that $xy$ is $\Integ$-homogeneous then both $x$ and $y$ are $\Integ$-homogeneous,
and the observation that $h-1$ is not $\Integ$-homogeneous).
Since $\rho(h)$ is a unit of $S$, $\rho$ factors through the localization $B \to B_h$. Then we have
$$
\xymatrix{
B \ar[r] \ar@/^1pc/[rr]^{\rho}  &  B_h \ar[r]_{\rho_h}  &  S \\
& B_{(h)} \ar@{^(->}[u] \ar[r]_{\rho'} & S_0 \ar@{^(->}[u]
}
$$
where (since $\rho(h)=1$) $\rho_h( x/h^n ) = \rho(x)$  for every $x \in B$ and $n\in\Nat$.
We claim that $\rho'$ is bijective.
Indeed, suppose that $x/h^n \in B_{(h)}$ ($n \in\Nat$, $x \in B_{mn}$) satisfies $\rho'(x/h^n) = 0$;
then $0 = \rho_h(x/h^n) = \rho(x)$, so $x =0$ by \eqref{okcjbv8238iwf}.
So $\rho'$ is injective.
To prove surjectivity, consider $y \in S_0$. 
Since $\rho( B_j' ) = S_j$ for each $j \in G$,  there exists $x \in B_0' = \bigoplus_{i \in m\Integ} B_i$ such that $\rho(x)=y$.
Write $x = \sum_{i \in \Integ} x_{mi}$ ($x_{mi} \in B_{mi}$);
then the element $\xi = \sum_{i \in \Integ} \frac{x_{mi}}{h^i}$ of $B_{(h)}$ satisfies
$\rho'(\xi) =  \sum_{i \in \Integ} \rho'\big( \frac{x_{mi}}{h^i} \big) =  \sum_{i \in \Integ} \rho_h\big( \frac{x_{mi}}{h^i} \big)
=  \sum_{i \in \Integ} \rho ( x_{mi} ) = \rho(x) = y$. 
So, 
\begin{equation} \label {kjbh9298whd0kspod}
\text{$\rho'$ is bijective and consequently $B_{(h)} \isom S_0$.}
\end{equation}

Since $h-1 \in A$, $D$ induces a $G$-homogeneous locally nilpotent derivation $\bar D : S \to S$; more precisely,
we have $\bar D( \rho(x) ) = \rho( D(x) )$ for all $x \in B$.
Let $R = \ker(\bar D)$ and note that $R$ is a $G$-graded subring of $S$.
By definition of $h$, there exists a $\Integ$-homogeneous $\beta \in B$ satisfying $D(\beta)=h$; define $v = \rho(\beta)$ and note that $v$ is a 
$G$-homogeneous element of $S$  satisfying $\bar D(v) = \bar D( \rho(\beta) ) = \rho( D(\beta) ) = \rho(h) = 1$.
Since $\Rat \subseteq B$ by assumption, the Slice Theorem \ref{SliceThm} implies that
$S = R[v] = R^{[1]}$. We record this as:
\begin{equation}  \label {do921fKF5Jjjanbu56yF6X}
\text{$S = R[v] = R^{[1]}$ for some $G$-homogeneous element $v$ of $S$.}
\end{equation}
Next, we claim that
\begin{equation}  \label {kcjnvop0293red}
\text{the set $H = \setspec{ i \in G }{ R_i \neq 0 }$ is a subgroup of $G$.}
\end{equation}
Note that if $R$ is a domain then $H$ is closed under addition and hence is a subgroup of $G$ (because $G$ is a finite group).
However, $R$ is not necessarily a domain.
The proof of \eqref{kcjnvop0293red} has several steps. First, we show that
\begin{equation}  \label {9w0902378687239rufri9e}
\text{if $g \in B$ satisfies $(h-1)g \in D(B)$, then $g \in D(B)$.}
\end{equation}
Consider $g \in B\setminus\{0\}$ such that $(h-1)g \in D(B)$.
Recall that $m \neq 0$.
If $m>0$ (resp.\ $m<0$), write $g = g_{i_0} + \sum_{i>i_0} g_i$ (resp.~$g = g_{i_0} + \sum_{i<i_0} g_i$),
where $g_i \in B_i$ and $g_{i_0} \in B_{i_0} \setminus \{0\}$.
Then $D(B) \ni (h-1)g = - g_{i_0} + w$ where $w = \sum_{i > i_0} w_i$ (resp.\ $w = \sum_{i < i_0} w_i$), $w_i \in B_i$.
So $g_{i_0} \in D(B)$ (because $D(B)$ is a {\it graded\/} $A$-submodule of $B$), which proves the claim in the special case where $g$ is $\Integ$-homogeneous.
If $g$ is not $\Integ$-homogeneous then $(h-1)g_{i_0} \in D(B)$ (because $g_{i_0} \in D(B)$ and $D(B)$ is an $A$-module),
so $(h-1)(g-g_{i_0}) \in D(B)$.
By induction on the number of nonzero $\Integ$-homogeneous terms in $g$, this proves \eqref{9w0902378687239rufri9e}.
Next, we show that
\begin{equation}  \label {eo9ivno9wds}
\rho(A) = R .
\end{equation}
Indeed, let $y \in R$.
Pick any $x \in B$ satisfying $\rho(x)=y$. Then $0 = \bar D(y) = \rho( D(x) )$, so $D(x) = (h-1)g$ for some $g \in B$.
Then $g \in D(B)$ by \eqref{9w0902378687239rufri9e}, so there exists $w \in B$ such that $D(w)=g$.
Since $D( x-(h-1)w ) = D(x) - (h-1)g=0$, we have $x-(h-1)w \in A$ and hence $y = \rho(x-(h-1)w) \in \rho(A)$.
Thus, $R \subseteq \rho(A)$.  It is clear that $\rho(A) \subseteq R$, so \eqref{eo9ivno9wds} follows.

To prove \eqref{kcjnvop0293red}, consider the submonoid $H' = \setspec{ i \in \Integ }{ A_i \neq 0 }$ of $\Integ$.
If $i \in H'$ then pick $a \in A_i \setminus \{0\}$ and note that $\rho(a)\neq 0$ by \eqref{okcjbv8238iwf}, so $\rho(a) \in R_{\pi(i)}\setminus\{0\}$,
 so $\pi(i) \in H$, which shows that $\pi(H') \subseteq H$.  If $j \in H$ then pick $y \in R_j \setminus \{0\}$;
by \eqref{eo9ivno9wds}, there exists $x \in A$ such that $\rho(x)=y$.
Write $x = \sum_{i \in \Integ} x_i$ with $x_i \in A_i$. Since $y\neq0$ and $\rho(B_i) \subseteq S_{\pi(i)}$ for every $i \in \Integ$, 
there exists an $i \in \Integ$ such that $\pi(i)=j$ and $x_i \in A_i\setminus\{0\}$; then $i \in H'$, so $j \in \pi(H')$,
which shows that $H \subseteq \pi(H')$.
Thus, $H  = \pi(H')$ is closed under addition. Since $G$ is finite, $H$ is a subgroup of $G$.

This proves \eqref{kcjnvop0293red}.
In view of \eqref{do921fKF5Jjjanbu56yF6X} and \eqref{kcjnvop0293red},
Lemma \ref{dpcinvbi293jhdbcnc38i} implies that one of the following conditions is satisfied:
\begin{enumerata}

\item $S_0 = (R_0)^{[1]}$,

\item there exists $e \in H \setminus \{0\}$ such that, if $k>0$ denotes the order of $e$ in $H$,
then, for every $r \in R_e \setminus\{0\}$, 
$(S_0)_{r^k}$ is a polynomial ring in one variable.

\end{enumerata}
Moreover, we have $B_{(h)} \isom S_0$ by \eqref{kjbh9298whd0kspod}.
If (a) holds then $B_{(h)}$ is a polynomial ring in one variable; since the $\Integ$-degree of $h$ is $m\neq0$, $h$ is a cylindrical element of $B$ and we are done.

Suppose that (b) holds, and let $e \in H \setminus \{0\}$ and $k>0$ be as in (b).  
Since $\pi(H')=H$, there exists $i \in H'$ such that $\pi(i)=e$.
We have $A_i \neq 0$ by definition of $H'$; pick $x \in A_i\setminus\{0\}$ and define $r = \rho(x) \in R_e \setminus\{0\}$
($r \neq 0$ by \eqref{okcjbv8238iwf}).
Let $a = x^k$; then $a$ is a nonzero $\Integ$-homogeneous element of $A$ and $\rho(a) = r^k$.
The $\Integ$-degree of $a$ is $ki \neq 0$, and $ki \in m\Integ$ because $\pi(ki)=k\pi(i)=ke=0$.
So there exists an integer $\ell>0$ such that some element $\alpha$ of $\{ a h^\ell, \frac{a}{h^\ell}\}$ belongs to $A_{(h)}$.
Since $\rho'(\alpha) = \rho(a)=r^k$, the isomorphism $\rho'$ extends to $(B_{(h)})_\alpha \isom (S_0)_{r^k}$, which, by (b), is a polynomial ring in one variable.
Since $B_{(ah)} \isom (B_{(h)})_\alpha$, $B_{(ah)}$ is a polynomial ring in one variable.
We have $B_{(ah)} = B_{(a^2h)}$ and at least one element $f$ of $\{ah, a^2h\}$ has nonzero $\Integ$-degree.
Then $f$ is a cylindrical element of $B$ and an element of $D(B) \cap \ker(D)$. This proves the Proposition.
\end{proof}

\begin{bigremark} \label {aPOKJxcbcdgfjkawEv3mnCigwvgdEduh9w}
We point out an issue with Lemma 1.8 of \cite{KPZ2013}.
That Lemma is used in the proof of Proposition 2.7, which is the crucial step in the proof of Theorem 0.6.
To prove the Lemma, it has to be shown that $F^{ \bar\partial } \subseteq \rho( A^{\partial} )$,
or equivalently $\rho^{-1}\big( F^{ \bar\partial } \big) \subseteq A^{\partial}$.
Let $E = \rho^{-1}\big( F^{ \bar\partial } \big)$, which is a subring of $A=\bigoplus_{j \in \Nat} A_j$.
The proof of  $E \subseteq A^{\partial}$ given in \cite{KPZ2013} consists in showing that $E \cap A_j \subseteq A^{\partial}$ for every $j$. 
For this to prove the desired inclusion, one needs
\begin{equation}  \label {lkjhbcvu63wte923u0}
\textstyle E = \sum_{j \in \Nat} (E \cap A_j) 
\end{equation}
to be true, but \eqref{lkjhbcvu63wte923u0} is not mentioned in the proof. 
In fact \eqref{lkjhbcvu63wte923u0} is true, but (as far as we can tell) is harder to prove than  $E \cap A_j \subseteq A^{\partial}$.
Indeed, $E \cap A_j \subseteq A^{\partial}$ follows from the simple observation \eqref{okcjbv8238iwf} in our proof of Proposition \ref{ckjb20932gediqw0idAwk},
whereas a proof of \eqref{lkjhbcvu63wte923u0} would have to include  something equivalent to our proof of \eqref{9w0902378687239rufri9e}.
\end{bigremark}

\begin{corollary}  \label {ekdh9c8wehwsgdiqolaj565v44rs7}
Let $G$ be an abelian group and $B = \bigoplus_{i \in G} B_i$ a $G$-graded $\Rat$-domain.
Suppose that $d$ is an element of $G$ of infinite order and that $D : B^{(d)} \to B^{(d)}$ is a nonzero homogeneous locally nilpotent derivation
such that $\ker D \nsubseteq B_0$.
Then some element of $D(B^{(d)}) \cap \ker(D)$ is a cylindrical element of $B$.
\end{corollary}

\begin{proof}
Since $\langle d \rangle \isom \Integ$,
we may regard $B^{(d)}  = \bigoplus_{i \in \langle d \rangle} B_i$ as a $\Integ$-graded ring.
Then Proposition \ref{ckjb20932gediqw0idAwk} implies that
some $f \in D(B^{(d)}) \cap \ker(D)$ is a cylindrical element of $B^{(d)}$.
Then $f$ is a nonzero homogeneous element of $B$ and $\deg(f)$ has infinite order.
Since $B_{(f)} = (B^{(d)})_{(f)}$, it follows that $B_{(f)}$ is a polynomial ring in one 
variable, so $f$ is a cylindrical element of $B$.
\end{proof}

\begin{proposition}  \label {do9uh298wcg8239eue}
Let $G$ be an abelian group and $B = \bigoplus_{i \in G} B_i$ a $G$-graded domain which is finitely generated 
as a $\bk$-algebra, where $\bk$ is a field of characteristic zero and $\bk \subseteq B_0$.
Let $f$ be a cylindrical element of $B$ and let $d = \deg(f) \in G$.
Then there exists a nonzero homogeneous locally nilpotent derivation $D : B^{(d)} \to B^{(d)}$ such
that $\ker(D) \nsubseteq B_0$ and, for some $n\ge1$, $f^n \in D(B^{(d)}) \cap \ker(D)$.
\end{proposition}

\begin{proof}
Since $B_{(f)}$ is a polynomial ring in one variable and is not the zero ring,
there exist a locally nilpotent derivation $\delta : B_{(f)} \to B_{(f)}$ and an element $t \in B_{(f)}$ such that $\delta(t)=1$
(write $B_{(f)} = R[t]=R^{[1]}$ and let $\delta : R[t]\to R[t]$ be the $t$-derivative).
The fact that $d = \deg(f)$ has infinite order implies that $B_{(f)}[f] = ( B_{(f)} )^{[1]}$,
so there is a unique locally nilpotent derivation $\Delta : B_{(f)}[f] \to B_{(f)}[f]$
that extends $\delta$ and satisfies $\Delta(f)=0$. By localization, we obtain 
a locally nilpotent derivation $\Delta' : B_{(f)}[f,f^{-1}] \to B_{(f)}[f,f^{-1}]$ that extends $\Delta$.
Note that $f \in B^{(d)}$ and that $(B^{(d)})_f = B_{(f)}[f,f^{-1}]$.
So $\Delta'$ is a locally nilpotent derivation of $(B^{(d)})_f$.
Since $B$ is a finitely generated $\bk$-algebra, so is $B^{(d)}$ by Lemma \ref{kcjviyc2e2w45w6d2o9kmgvczxwer}. 
As is well known, this implies that there exists a positive integer $n$ such that $f^n \Delta' : (B^{(d)})_f \to (B^{(d)})_f$
maps $B^{(d)}$ into itself.
Since $f^n\Delta'$ is locally nilpotent (because $\Delta'(f)=0$), it follows that
the restriction $D :  B^{(d)} \to B^{(d)}$ of $f^n\Delta'$ is a locally nilpotent derivation.
Moreover, it is easy to see that $D$ is homogeneous. Since $D(f)=0$, we have $\ker(D) \nsubseteq B_0$.
Since $t \in B_{(f)}$, we can pick $m\ge1$ such that $tf^m \in B^{(d)}$.
Then $D(tf^m) = f^n\Delta'(tf^m) = f^{m+n} \delta(t) = f^{m+n}$, so $D\neq 0$ and $f^{m+n} \in D(B^{(d)}) \cap \ker(D)$.
\end{proof}


\begin{nothing*} \label {ajsdfo93ew0sdpw0edp}
{\bf Proof of Theorem \ref{edh83yf6r79hvujhxu6wrefji9e}. }
Let us prove assertion (1).
Let $\bk$ be the field over which $B$ is a finitely generated algebra.
Since $\bk \subseteq B$ and $B$ is a domain graded by a torsion-free abelian group, we have $\bk \subseteq B_0$.
So (b) implies (a) by Proposition \ref{do9uh298wcg8239eue}.
We prove the converse.
In view of Corollary \ref{ekdh9c8wehwsgdiqolaj565v44rs7},
it suffices to show that
\begin{equation} \label {lkjbo928dA0qbvd7q23j}
\begin{minipage}[t]{.85\textwidth}
if $d\in \Integ \setminus\{0\}$ and $B^{(d)}$ is non-rigid, then
there exists a nonzero homogeneous locally nilpotent derivation $D : B^{(d)} \to B^{(d)}$ such that $\ker(D) \nsubseteq B_0$.
\end{minipage}
\end{equation}

Suppose that $d \neq0$ is such that  $B^{(d)}$ is non-rigid.
Note that $B^{(d)}$ is a $\Integ$-graded domain which (by Lemma \ref{kcjviyc2e2w45w6d2o9kmgvczxwer}) is a finitely generated $\bk$-algebra.
By Lemma \ref{dkjcnbv293qeokwdncbow90}, there exists a nonzero homogeneous locally nilpotent derivation
$D : B^{(d)} \to B^{(d)}$. We have to show that $\ker(D) \nsubseteq B_0$.
Arguing by contradiction, assume that  $\ker(D) \subseteq B_0$.
Since (by Lemma \ref{dkjhf2983ed9eje9}) $\trdeg(B : \ker D ) = 1$, condition (i) cannot hold; so (ii) holds and we may consider
$x \in B_i\setminus\{0\}$ and $y \in B_j\setminus\{0\}$ where $i < 0 < j$.
Clearly, $y$ is transcendental over $B_0$; so $\trdeg(B : B_0) \ge 1 = \trdeg(B: \ker D)$,
so $B_0$ is algebraic over $\ker(D)$, so $B_0 = \ker(D)$, since $\ker(D)$ is algebraically closed in $B$.
We have $x^jy^{|i|} \in B_0 = \ker(D)$, so $x,y\in\ker(D)$ because $\ker(D)$ is factorially closed in $B$ by Lemma \ref{dkjhf2983ed9eje9}.
So $x,y \in B_0$, which contradicts the choice of $x,y$.
So $\ker(D) \nsubseteq B_0$.
This proves \eqref{lkjbo928dA0qbvd7q23j} and completes the proof of assertion (1) of the Theorem.
Assertion (2) immediately follows from assertion (1) and Theorem \ref{kjchvtdes334eOF5A5WF6I4ut9ity0o6plkg}
(or its special case Corollary \ref{eocnb9nxrd4g5sywjaafgbb}).  \hfill\qedsymbol
\end{nothing*}

We use the following notation in the next proof.
If $B$ is a $\Integ$-graded domain, we write $\HFrac(B)$ for the degree-$0$ subring of the $\Integ$-graded ring $S^{-1}B$,
where $S$ is the set of all nonzero homogeneous elements of $B$.
Note that $\HFrac(B)$ is a subfield of $\Frac B$, and that if the grading is non-trivial then $\Frac(B) = \big( \HFrac B \big)^{(1)}$.
Also, it is easy to check that $\HFrac(B^{(d)}) = \HFrac(B)$ for any $d \in \Integ\setminus\{0\}$.

\begin{nothing*} \label {cjOlKo10eof9ccnqpkeuycazvds4wud}
{\bf Proof of Proposition \ref{kcnboiqwdcm203}. }
Let $f$ be a cylindrical element of $B$.  Since $f$ is a nonzero homogeneous element of nonzero degree, the grading is non-trivial and
consequently $\Frac(B) = \big( \HFrac B \big)^{(1)}$.
Let $d = \deg(f)$.
Since $\Frac( B_{(f)} ) = \HFrac(B^{(d)}) = \HFrac(B)$, we have
$$
\Frac(B) = \big( \Frac B_{(f)} \big)^{(1)} .
$$
Since $B_{(f)} = R^{[1]}$ for some subring $R$ of $B_{(f)}$, we have $\Frac(B) = K^{(2)}$ where $K = \Frac(R)$.
Since $B$ is a domain graded by a torsion-free abelian group, we have $\bk \subseteq B_0$. So $\bk \subseteq B_{(f)} = R^{[1]}$,
so $\bk \subseteq R$, so $\bk \subseteq K$.  \hfill\qedsymbol
\end{nothing*}


\section{Cylindrical elements and $H$-polar cylinders}
\label {jhfbo293wei0dc02uew9}

The purpose of this section is to show that Theorem 0.6 and Corollary 3.2 of \cite{KPZ2013} are equivalent to assertion \ref{iocjvn2uygew8dhc},
stated in the introduction.
We begin by recalling the part of the Dolgachev-Pinkham-Demazure (DPD) construction that we need to properly state the results of \cite{KPZ2013}.

\begin{nothing*} \label {kcmvo2938wndsco09}
Let $Y$ be a noetherian normal integral scheme and $K$ the function field of $Y$. 
We write $\Div(Y)$ for the group of Weil divisors of $Y$.
A divisor $D \in \Div(Y)$ is said to be {\it Cartier\/} if $Y$ can be covered by open sets $U_i$ such that, for each $i$,
we have $D|_{U_i} = \div_{U_i}(f_i)$ for some $f_i \in K^*$.
If $D \in \Div(Y)$ then the sheaf $\Oeul_Y(D)$ on $Y$ is defined by
$$
\Gamma(U,\Oeul_Y(D)) = \{0\} \cup \setspec{ f \in K^* }{ \div_U(f) + D|_U \ge 0 } \qquad \text{($\emptyset \neq U \subseteq Y$ open)}.
$$
It is well known that $\Oeul_Y(D)$ is coherent, and that $D$ is Cartier if and only if $\Oeul(D)$ is invertible.
Refer to \cite[II.7.4]{Hartshorne} for the notion of an {\it ample invertible sheaf\/} on $Y$.
One says that the divisor $D \in \Div(Y)$ is {\it ample\/} if $\Oeul(D)$ is an ample invertible sheaf.
It follows that ample divisors are in particular Cartier.
\end{nothing*}

\begin{nothing*}
Given $Y$ and $K$ as in \ref{kcmvo2938wndsco09},
let $\Div_\Rat(Y)$ be the group of $\Rat$-divisors of $Y$, i.e., the free $\Rat$-module on the set of prime divisors of $Y$
(so $\Div(Y) \subseteq \Div_\Rat(Y)$).
Given $D = \sum_i r_i C_i \in \Div_\Rat(Y)$ (where the $C_i$ are prime divisors and $r_i \in \Rat$), one defines
$\lfloor D \rfloor = \sum_i \lfloor r_i \rfloor C_i \in \Div(Y)$.
Elements $D$ and $D'$ of $\Div_\Rat(Y)$ are said to be {\it linearly equivalent\/} ($D \sim D'$)
if there exists  $f \in K^*$ such that $D - D' = \div_Y(f)$.\footnote{Different authors use different definitions of linear equivalence in $\Div_\Rat(Y)$.
We use the same definition as in \cite{KPZ2013} and \cite{Demazure_1988}.}
A $\Rat$-divisor $D \in \Div_\Rat(Y)$ is said to be {\it $\Rat$-ample\/} if there exists $n\ge1$ such that $nD \in \Div(Y)$
and $nD$ is ample in the sense of \ref{kcmvo2938wndsco09}.
Note that a divisor $D \in \Div(Y)$ is ample if and only if it is Cartier and $\Rat$-ample.
\end{nothing*}

\begin{nothing*} \label {pq09i34nb0q39efi023}
Let $B$ be an $\Nat$-graded domain.
An element $\xi$ of $\Frac B$ is said to be {\it homogeneous\/} if it can be written as $\xi = a/b$ for some homogeneous elements $a,b \in B$ with $b \neq 0$.
Moreover, if $\xi=a/b$ with $a \in B_m$ and  $b \in B_n \setminus \{0\}$ then we say that $\xi$ is homogeneous of degree $m-n$. 
We write $\big( \Frac B \big)_d$ for the set of homogeneous elements of $\Frac B$ of degree $d$.
Note that the function field of $\Proj B$ is $\big(\Frac B\big)_0$. 
We shall also use the following notation, for a homogeneous element $f$ of $B$:
$$
\bbV_+(f) = \setspec{ \pgoth \in \Proj(B) }{ f \in \pgoth } \quad \text{and} \quad \bbD_+(f) = \Proj(B) \setminus \bbV_+(f) .
$$
\end{nothing*}

\begin{nothing*} \label {Pkcnbvc9w3eidjojf0q9w}
Let $B$ be an $\Nat$-graded noetherian normal domain such that the prime ideal $B_+$ has height greater than $1$.
Let $X = \Spec B$ and $Y = \Proj B$.
We shall now define a $\Rat$-linear map $D \mapsto D^*$ from $\Div_\Rat(Y)$ to $\Div_\Rat(X)$.

Let $K(X)$ and $K(Y)$ be the function fields of $X$ and $Y$ respectively.
Let $\YY$ be the set of homogeneous prime ideals of $B$ of height $1$.
Since $\haut(B_+)>1$, we have $\YY = \setspec{ y \in Y }{ \dim\Oeul_{Y,y}=1 }$.
For each $\pgoth \in \YY$, $B_\pgoth \supset B_{(\pgoth)}$ is an extension of DVRs; let $e_\pgoth$ denote the ramification index of this extension.
Then $e_\pgoth \in \Nat\setminus \{0\}$.
If $v^Y_\pgoth : K(Y)^* \to \Integ$ and 
$v^X_\pgoth : K(X)^* \to \Integ$ denote the normalized\footnote{The word ``normalized'' means that the maps $v^Y_\pgoth$ and $v^X_\pgoth$ are surjective.}
valuations of $B_{(\pgoth)}$ and $B_\pgoth$ respectively,
then $v^X_\pgoth (\xi) = e_\pgoth v^Y_\pgoth(\xi)$ for all $\xi \in K(Y)^*$.
Let $C_\pgoth^Y$ (resp.~$C_\pgoth^X$) denote the closure of $\{ \pgoth \}$ in $Y$ (resp.\ in $X$).
Then $C_\pgoth^Y$ (resp.~$C_\pgoth^X$) is a prime divisor of $Y$ (resp.\ of $X$),
and every prime divisor of $Y$ is a $C^Y_\pgoth$ for some $\pgoth \in \YY$.
We define $( C_\pgoth^Y )^* = e_\pgoth C_\pgoth^X$ for each $\pgoth \in \YY$,
and extend linearly to  a $\Rat$-linear map $\Div_\Rat(Y) \to \Div_\Rat(X)$, $D \mapsto D^*$.
It is not hard to see that the linear map $D \mapsto D^*$ has the following two properties:
\begin{gather}
\label {PWSoidbfci27wyd}
\begin{minipage}[t]{.85\textwidth}
\it  $\big( \div_Y(\xi) \big)^* = \div_X(\xi)$ for all $\xi \in K(Y)^*$,
\end{minipage} \\
\label {Dpojfnb239wdwkd}
\begin{minipage}[t]{.85\textwidth}
\it if $f$ is a nonzero homogeneous element of $B$ and $D \in \Div_\Rat(Y)$ satisfies
$D^* = \div_X(f)$, then $D \ge 0$ and $\supp(D) = \bbV_+(f)$.
\end{minipage}
\end{gather}
\end{nothing*}

\begin{notation} \label {ed9e98e7326531527ydfjn8rn73}
Given a noetherian normal integral scheme $Y$ and $D \in \Div_\Rat(Y)$, define the $\Nat$-graded ring
$$
A(Y,D) = \bigoplus_{i \in \Nat} H^0(Y,\Oeul_Y(\lfloor iD \rfloor)) T^i ,
$$
where $T$ is an indeterminate over the function field $K$ of $Y$. Note that $A(Y,D)$ is a graded subring of $K[T]$ and that $A_0 = \Oeul_Y(Y)$.
\end{notation}

\begin{notation} \label {ed988tggfweCYHbj847r5823dgw89}
Given a field $\bk$, let $\PP$ be the class of pairs $(Y,H)$ satisfying (i) and (ii):
\begin{itemize}

\item[(i)] $Y$ is a normal integral scheme which is projective over a finitely generated $\bk$-algebra $R$ such that $\dim Y > \dim R$;

\item[(ii)] $H$ belongs to $\Div_\Rat(Y)$ and is $\Rat$-ample.

\end{itemize}
Note that (i) has the following consequence:
\begin{itemize}

\item[(iii)] $Y$ is of finite type over $\bk$, $\Oeul_Y(Y)$ is a finitely generated $\bk$-algebra such that $\dim Y > \dim\Oeul_Y(Y)$,
and the canonical morphism $Y \to \Spec \Oeul_Y(Y)$ is projective.

\end{itemize}
(By \cite[Thm III.5.2]{Hartshorne}, $\Oeul_Y(Y)$ is a finitely generated $R$-module;
so $\Oeul_Y(Y)$ is a finitely generated $\bk$-algebra and $\dim Y > \dim R \ge \dim\Oeul_Y(Y)$.
Since the composition $Y \to \Spec\Oeul_Y(Y) \to \Spec R$ is projective and $\Spec\Oeul_Y(Y) \to \Spec R$ is separated,
$Y \to \Spec \Oeul_Y(Y)$ is projective and (iii) is true.)
\end{notation}

The following is well known.

\begin{theorem}  \label {kjxbco23idjq0wds0}
Let $(Y,H) \in \PP$ and $A=A(Y,H)$.
\begin{enumerata}

\item The ring $A$ is an $\Nat$-graded normal domain, a finitely generated $\bk$-algebra, and satisfies
$\Frac( A ) = K(T)$, $e( A ) = 1$, and $\haut(A_+)>1$, where $K$ is the function field of $Y$.

\item Let $Z = \Proj A$ and $X = \Spec A$, and consider the map
$\Div_\Rat(Z) \to \Div_\Rat(X)$, $D \mapsto D^*$, defined in \ref{Pkcnbvc9w3eidjojf0q9w}.
There exists an isomorphism $j : Y \to Z$ of schemes over $\bk$ such that,
if $D \in \Div_\Rat(Z)$ is the image of $H$ by the isomorphism $\Div_\Rat(Y) \to \Div_\Rat(Z)$ induced by $j$, then $D^* = \div_X(T)$.

\item If $H' \in \Div_\Rat(Y)$ is linearly equivalent to $H$ then $A(Y,H)$ and $A(Y,H')$ are isomorphic as graded $\bk$-algebras.

\end{enumerata}
\end{theorem}

\begin{nothing*}
Here are some references for Theorem \ref{kjxbco23idjq0wds0}.
(a) The Corollaire in \cite[3.2]{Demazure_1988} shows that $\Frac( A ) = K(T)$, which implies that  $e( A ) = 1$.
The fact that $A$ is normal is stated without proof in \cite[3.1]{Demazure_1988}; the proof can be found in that of the Lemme in \cite[2.7]{Demazure_1988}.
The fact that $A$ is a finitely generated $\bk$-algebra is proved in Proposition 3.3 of \cite{Demazure_1988},
apparently under the assumption that $Y$ is projective over $\bk$. The general case is well known, and can be seen as follows.
Recall that $Y$ is projective over a finitely generated $\bk$-algebra $R$.
Let $N>0$ be such that $NH \in \Div(Y)$ is an ample Cartier divisor.
Let $\Leul = \Oeul_Y( \lfloor NH \rfloor ) = \Oeul_Y( NH )$.
Then $\Leul$ is an ample invertible sheaf on $Y$, so
$\bigoplus_{n \in \Nat} \Gamma(Y, \Leul^{\otimes n})$ is a finitely generated $R$-algebra (hence a finitely generated $\bk$-algebra)
by Proposition 3.2.1 of \cite{Dolg:McKayCorr}.
Since $A^{(N)} = \bigoplus_{n \in \Nat} \Gamma(Y, \Leul^{\otimes n})T^{Nn}$, $A^{(N)}$ is a finitely generated $\bk$-algebra.
The fact that $A$ is a finitely generated $\bk$-algebra then follows by the argument given in the proof of Proposition 3.3 of \cite{Demazure_1988}.
Finally, we have $\haut(A_+) = \dim A - \dim A_0 = (\dim(Y)+1) - \dim \Oeul_Y(Y) > 1$, where we use
$A_0 = \Oeul_Y(Y)$, $\dim Y > \dim \Oeul_Y(Y)$ (see \ref{ed988tggfweCYHbj847r5823dgw89}(iii)), and $\dim(A) = \dim(Y)+1$ (because $Y \isom \Proj A$, see (b)). 

(b) The Proposition in paragraph 3.2 of \cite{Demazure_1988} defines an open immersion $j : Y \to Z$ which is such that
$$
\xymatrix@R=15pt{
Y \ar@{^{(}->}[r]^-{j} \ar[d]_-{p_Y} & Z  \ar[d]^-{p_Z} \\
\Spec\Oeul_Y(Y) \ar@{=}[r] &  \Spec A_0
}
$$
commutes.
It was noted in \ref{ed988tggfweCYHbj847r5823dgw89}(iii) that the canonical morphism $p_Y$ is projective;
so $p_Z \circ j$ is projective, hence proper. Since $p_Z$ is separated, $j$ is proper, so $j$ is an isomorphism.
Commutativity of the diagram also shows that $j$ is a morphism over $\bk$, since $\bk \subseteq A_0$.
The equality $D^* = \div_X(T)$ follows from the Corollaire in \cite[2.9]{Demazure_1988}, together with the
remark on page 50 that the results of 2.9 ``se transportent au c\^one \underline{mutatis mutandis}''.

(c) If $H'=H+\div_Y(\xi)$, where $\xi \in K^*$, then the $K$-automorphism of $K[T]$ that sends $T$ to $\xi T$ 
maps $H^0(Y,\Oeul_Y(nH'))T^n$ onto $H^0(Y,\Oeul_Y(nH))T^n$ (for each $n\in\Nat$).
Since $\bk \subseteq K$, this gives a $\bk$-isomorphism of graded rings $A(Y,H') \to A(Y,H)$.
\end{nothing*}

\begin{theorem}[Th\'eor\`eme 3.5, \cite{Demazure_1988}] \label {skjbci8wjqpodmpa0}
Let $\bk$ be a field and $B$ an $\Nat$-graded normal domain which is a finitely generated $\bk$-algebra and such that $e(B)=1$ and $\haut(B_+)>1$.
Let $u \in ( \Frac B )_1 \setminus \{0\}$.
There exists a unique $H \in \Div_\Rat(Y)$ satisfying $H^* = \div_X(u)$, where $Y = \Proj B$ and $X = \Spec B$.
Moreover, $H$ is $\Rat$-ample and 
$$
B = \bigoplus_{n \in \Nat} H^0( Y, \Oeul( \lfloor n H \rfloor ) ) u^n.
$$
\end{theorem}

\begin{bigremarks}  \label {jcbvio8q2wdxjvyQ6}
The following comments are related to Theorem \ref{skjbci8wjqpodmpa0}.
\begin{enumerata}

\item The assumption $e(B)=1$ implies that $( \Frac B )_1 \setminus \{0\} \neq \emptyset$.

\item The last assertion of the Theorem is an {\it equality\/} of rings.

\item Since $B$ is a finitely generated $\bk$-algebra, so is $B_0 \isom B/B_+$.
The condition $\haut(B_+)>1$ implies that $\dim Y > \dim B_0$.
Since the canonical morphism $\Proj B \to \Spec B_0$ is projective, we see that $Y$ is projective over a finitely generated $\bk$-algebra $R=B_0$
such that $\dim Y>\dim R$.  That is, $(Y,H) \in \PP$. Moreover, $A(Y,H) \isom B$.

\item If $u,u' \in ( \Frac B )_1 \setminus \{0\}$ and $H,H' \in \Div_\Rat(Y)$ satisfy $H^* = \div_X(u)$ and $(H')^* = \div_X(u')$,
then $H \sim H'$. Indeed, $H-H' = \div_Y(u/u')$.

\end{enumerata}
\end{bigremarks}

\begin{nothing*} \label {09203ndc020B0ej102}
Let $\bk$ be a field.
Theorems \ref{kjxbco23idjq0wds0} and \ref{skjbci8wjqpodmpa0} define two bijections, inverse of each other:
\begin{enumerate}

\item Define an equivalence relation $\approx$ on $\PP$ by declaring that  $(Y,H) \approx (Y',H')$ if and only if
there exists an isomorphism $\theta : Y \to Y'$ over $\Spec \bk$ that carries $H$ to a $\Rat$-divisor of $Y'$ linearly equivalent to $H'$.
The equivalence class of $(Y,H)$ is denoted $[Y,H]$.

\item Let $\GND$ be the class of $\Nat$-graded normal domains $B$
such that $e(B)=1$, $\haut(B_+)>1$, and $B$ is a finitely generated $\bk$-algebra.
Two elements $B = \bigoplus_{i\in\Nat} B_i$ and $B' = \bigoplus_{i\in\Nat} B_i'$
of $\GND$ are {\it isomorphic\/} ($B \isom B'$) if there exists an isomorphism of $\bk$-algebras $\psi : B \to B'$ such that $\psi(B_i)=B_i'$ for all $i \in \Nat$.

\item Theorems \ref{kjxbco23idjq0wds0} and \ref{skjbci8wjqpodmpa0} and  Remark~\ref{jcbvio8q2wdxjvyQ6} give two maps
$$
\begin{array}[t]{rcl}
{\PP/\!\!\approx} & \xrightarrow{\text{\ \ref{kjxbco23idjq0wds0}\ }} & {\GND/\!\!\isom} \\
{[Y,H]} & \mapsto  & {[ A(Y,H) ]}
\end{array}
\quad\text{and}\quad
\begin{array}[t]{rcl}
{\GND/\!\!\isom} & \xrightarrow{\text{\ \ref{skjbci8wjqpodmpa0}\ }} & {\PP/\!\!\approx} \\
{[B]} & \mapsto & {[Y,H]} .
\end{array}
$$
By part (b) of Theorem \ref{kjxbco23idjq0wds0}, the composition
$$
\PP/\!\!\approx \,\,\, \longrightarrow \GND/\!\!\isom \,\,\, \longrightarrow \PP/\!\!\approx
$$
is the identity map of $\PP/\!\!\approx$.  By Theorem \ref{skjbci8wjqpodmpa0}, the composition
$$
\GND/\!\!\isom \,\,\, \longrightarrow \PP/\!\!\approx \,\,\, \longrightarrow \GND/\!\!\isom
$$
is the identity map of $\GND/\!\!\isom$.

\end{enumerate}
\end{nothing*}

\begin{caution} \label {IUpSIOdjbo93wew}
Let $\bk$ be a field, $B \in \GND$, $X = \Spec B$ and $Y = \Proj B$. Consider:
\begin{enumerata}

\item the set of  $\Rat$-ample $\Rat$-divisors $H$ of $Y$ such that $A(Y,H) \isom B$,
\item the set of $H \in \Div_\Rat(Y)$ such that $H^* = \div_X(u)$ for some $u \in \big( \Frac B \big)_1 \setminus \{0\}$.

\end{enumerata}
By Theorem \ref{skjbci8wjqpodmpa0}, (b) is a nonempty subset of (a).
By Remark \ref{jcbvio8q2wdxjvyQ6}(d), any two elements of (b) are linearly equivalent.
However, the elements of (a) are not necessarily linearly equivalent to one another (example: if $P,Q$ are distinct
closed points of $\proj^1$ then $A(\proj^1, \frac12 P ) \isom A(\proj^1, \frac12 Q )$ but $\frac12 P \nsim\frac12 Q$).
Thus, in some cases, (b) is a proper subset of (a).
However, injectivity of the map ${\PP/\!\!\approx} \to {\GND/\!\!\isom}$\ of \ref{09203ndc020B0ej102} implies that
if $H,H'$ are elements of (a) then there exists a $\bk$-automorphism of $Y$ that carries $H$ to a $\Rat$-divisor of $Y$ that is linearly equivalent to $H'$.
\end{caution}

Let us now define some terminology used in \cite{KPZ2013}.

\begin{definition} \label {pcoijvbo39urdmn38e}
Let $\bk$ be a field and $(Y,H) \in \PP$.
An open subset $U$ of $Y$ is {\it $H$-polar\/} if $U = Y \setminus \supp(D)$
for some effective $\Rat$-divisor $D \in \Div_\Rat(Y)$ which is linearly equivalent to $sH$ for some $s \in \Rat_{>0}$.
\end{definition}

\begin{bigremark} \label {c09fuhb2w90dsd}
Let  $(Y,H) \in \PP$ and let $H' \in \Div_\Rat(Y)$ be such that $\alpha H \sim \beta H'$ for some $\alpha, \beta \in \Rat_{>0}$.
Then $(Y,H') \in \PP$ and, for any open subset $U$ of $Y$,
$$
\text{$U$ is $H$-polar if and only if $U$ is $H'$-polar.}
$$
\end{bigremark}

\begin{definition} \label {ps9d8uu6ejdr3trb8}
Let $\bk$ be a field and $(Y,H) \in \PP$.
An open subset $U$ of $Y$ is a {\it cylinder\/} if $U \isom \aff^1 \times Z$ for some variety $Z$.
If $U$ is a cylinder and is $H$-polar in the sense of \ref{pcoijvbo39urdmn38e}, we call it an {\it $H$-polar cylinder}.
The pair $(Y,H)$ is {\it cylindrical\/} if some open subset of $Y$ is an $H$-polar cylinder.
\end{definition}

\begin{bigremark}  \label {xjkh9182djADBonqws7}
Let $(Y,H), (Y',H') \in \PP$ be such that $(Y,H) \approx (Y',H')$ (see \ref{09203ndc020B0ej102}).
Then $(Y,H)$ is cylindrical if and only if $(Y',H')$ is cylindrical.
(This follows from Remark \ref{c09fuhb2w90dsd}.)
\end{bigremark}

\begin{definition} \label {cvbo23ewpdw0k}
In the terminology of \cite{KPZ2013},
an affine $\bk$-variety $\Spec R$ is {\it cylindrical\/} if there exists $f \in R \setminus \{0\}$ such that 
$\bbD(f) \isom \aff_\bk^1 \times Z$ for some affine variety~$Z$.
\end{definition}

In view of the second paragraph of the Introduction, we see that $\Spec R$ is cylindrical if and only if $R$ is non-rigid (assuming that $\Char\bk=0$).
So, in Theorem 0.6, the condition ``$V$ is cylindrical'' is equivalent to $A$ being non-rigid,
and ``$V^{(d)}=\Spec A^{(d)}$ is cylindrical'' is equivalent to $A^{(d)}$ being non-rigid.

\medskip

\noindent{\bf Theorem 0.6 of \cite{KPZ2013}. }
[{\it Let $\bk$ be an algebraically closed field of characteristic $0$.}]
{\it
Let $A=\bigoplus_{\nu\ge 0} A_\nu$
be a positively graded affine domain over $\bk$.
Define the projective  variety $Y=\Proj A$
relative to this grading, let $H$  be the associated
$\Rat$-divisor on $Y$, and let $V=\Spec A$,
the affine quasicone over $Y$.
\begin{enumerata}

\item[(a)] If $V$ is cylindrical, then the associated pair $(Y,H)$ is  cylindrical.

\item[(b)] If the pair $(Y,H)$ is cylindrical, then for some $d\in\Nat$ the Veronese cone
$V^{(d)}=\Spec A^{(d)}$ is cylindrical, where $A^{(d)}=\bigoplus_{\nu\ge 0} A_{d\nu}$. 

\end{enumerata}
}

\medskip

\begin{nothing*}  \label {knxcbufk3idbcez55}
Three assumptions are missing from the statement of Theorem~0.6: (i)~$A$ is normal; (ii)~$e(A)=1$; and
(iii)~$A$ has transcendence degree at least $2$ over $A_0$, or equivalently, $\haut(A_+) > 1$.

Indeed, assumptions (i) and (ii) are used in the proof, and hence should appear in the statement of the Theorem.
Moreover, $A(Y,H)$ is always normal and such that $e\big( A(Y,H) \big) = 1$,
so if one of (i), (ii) is false then no $H$ can satisfy $A(Y,H)=A$, so the associated $\Rat$-divisor $H$ does not exist.

Here is an example showing that assertion (a) of Theorem 0.6 is false if we do not assume (iii).
Let $R$ be a normal affine $\bk$-domain such that no open subset of $\Spec R$ is isomorphic to a product $\aff^1 \times Z$ with $Z$ a variety.
Let $A = R[T] = R^{[1]}$, with standard $\Nat$-grading ($A_0=R$ and $T \in A_1$).
It is clear that $V = \Spec A$ is cylindrical.  Let $Y = \Proj A$.
Then $Y = \bbD_+(T) \isom \Spec A_{(T)} \isom \Spec R$, so $Y$ does not contain any open set of the form $\aff^1 \times Z$.
In particular, no element $H$ of $\Div_\Rat(Y)$ is such that $(Y,H)$ is cylindrical, so assertion (a) is false.
(In the proof of the Theorem, just before \cite[2.8]{KPZ2013}, the authors write that if (iii) is false then assertion (a) of Theorem 0.6
is true by virtue of Proposition 0.5.  By the above example, that is not correct.)
\end{nothing*}

\begin{nothing*}
We shall now restate Theorem 0.6 of \cite{KPZ2013} with the correct assumptions, i.e.,
we shall assume that $A$ satisfies conditions (i--iii) of paragraph \ref{knxcbufk3idbcez55}.
This is equivalent to assuming that $A \in \GND$ (see \ref{09203ndc020B0ej102}).
There remains the question of how to interpret the phrase ``let $H$ be the associated $\Rat$-divisor on $Y$'' (where $Y=\Proj A$).
A priori, it could mean that $H$ belongs to one or the other of the sets (a) and (b) of \ref{IUpSIOdjbo93wew}.
It turns out that the Theorem is valid regardless of which interpretation we choose, so we opt for the more general case
(i.e., $H$ belongs to the larger set (a)):
\end{nothing*}

\noindent{\bf Restatement of Theorem 0.6.}
{\it Let $\bk$ be an algebraically closed field of characteristic zero, let $A \in \GND$, 
let $Y = \Proj A$ and let $H$ be a $\Rat$-ample $\Rat$-divisor of $Y$ such that $A(Y,H) \isom A$.
Let $V=\Spec A$.
\begin{enumerata}

\item If $V$ is cylindrical, then $(Y,H)$ is  cylindrical.

\item If the pair $(Y,H)$ is cylindrical, then for some $d\in\Nat$ the Veronese cone $V^{(d)}=\Spec A^{(d)}$ is cylindrical.

\end{enumerata}}

\medskip

We want to understand what the condition ``$(Y,H)$ is cylindrical'' means for the graded ring $A$ in the above statement.
This is achieved in Corollary \ref{kjxcnopo93ewjdspwqW0}(a).

\begin{lemma}  \label {oijx9w38h8f9wlodl}
Let $\bk$ be a field, $B \in \GND$, $X = \Spec B$ and $Y = \Proj B$.
Let $u \in \big( \Frac B \big)_1 \setminus \{0\}$ and let $H \in \Div_\Rat(Y)$ be such that $H^* = \div_X(u)$.
\begin{enumerata}

\item An open subset $U$ of $Y$ is $H$-polar if and only if there exist $n\ge1$ and $f \in B_n \setminus\{0\}$ such 
that  $U = \bbD_+(f)$.

\item An open subset $U$ of $Y$ is an $H$-polar cylinder if and only if $U = \bbD_+(f)$ for some cylindrical element $f$ of $B$.

\item The pair $(Y,H)$ is cylindrical if and only if $B$ has a cylindrical element.

\item The following are equivalent.
\begin{enumerata}

\item $H$ is a Cartier divisor of $Y$;
\item $e(B/\pgoth) = 1$ for all $\pgoth \in \Proj B$.

\end{enumerata}

\end{enumerata}
\end{lemma}

\begin{proof} 
Theorem \ref{skjbci8wjqpodmpa0} implies that $u$ and $H$ exist and that we have the equality of rings
\begin{equation}  \label {98987291t7ydcsgd71}
\textstyle    B = \bigoplus_{n \in \Nat} H^0( Y, \Oeul( \lfloor n H \rfloor ) ) u^n.
\end{equation}
Let us first prove:
\begin{equation} \label {eicno9ewnd9Awojn}
\begin{minipage}[t]{.85\textwidth}
\it If $n \in \Nat \setminus \{0\}$ and $f \in B_n \setminus \{0\}$, then the $\Rat$-divisor $D = \div_Y(f/u^n) + nH$
satisfies $D\ge0$ and $\supp(D) = \bbV_+(f)$.
\end{minipage}
\end{equation}
To see this, first note that $f/u^n \in \big( \Frac B \big)_0 = K(Y)$, so it makes sense to define $D = \div_Y(f/u^n) + nH \in \Div_\Rat(Y)$.
We have $D^* = \big(\div_Y(f/u^n)\big)^* + nH^* = \div_X(f/u^n) + n \div_X(u) = \div_X(f)$ by \eqref{PWSoidbfci27wyd} and because $H^*=\div_X(u)$,
so \eqref{Dpojfnb239wdwkd} gives $D\ge0$ and $\supp(D) = \bbV_+(f)$, proving \eqref{eicno9ewnd9Awojn}.

(a) Let $U$ be an $H$-polar open set.
Then $U = Y \setminus \supp(D)$ for some $D \in\Div_\Rat(Y)$ such that $D\ge0$ and $D \sim sH$ for some $s \in \Rat_{>0}$.
Write $\div_Y(g) + sH = D \ge 0$, where $g \in K(Y)^*$.
Choose $m\in\Integ_{>0}$ such that $ms \in \Integ_{>0}$ and write $n = ms$.
Then $\div_Y(g^m) + nH = mD \ge 0$, so $\div_Y(g^m) + \lfloor nH \rfloor \ge 0$, so $g^m \in H^0(Y,\Oeul( \lfloor nH \rfloor ))$.
Let $f=g^mu^n$. We have $f \in B_n \setminus \{0\}$ by \eqref{98987291t7ydcsgd71}, and \eqref{eicno9ewnd9Awojn} then implies
that $\bbV_+(f) = \supp(  \div_Y( f u^{-n} ) + n H  ) = \supp(mD) = \supp(D)$,
so $U = Y \setminus \bbV_+(f) = \bbD_+(f)$.

Conversely, let  $n \in \Nat \setminus \{0\}$, $f \in B_n \setminus \{0\}$ and $U = \bbD_+(f)$.
By \eqref{eicno9ewnd9Awojn}, the $\Rat$-divisor $D = \div_Y(f/u^n) + nH$ satisfies
$D\ge0$ and $\supp(D) = \bbV_+(f)$, so $U=Y \setminus \supp(D)$.
Since $D$ is effective and linearly equivalent to $nH$, $U$ is $H$-polar. This proves (a).

(b) Suppose that $U$ is an $H$-polar cylinder. By part (a), we have $U = \bbD_+(f)$ for some nonzero homogeneous element $f$
of $B$ of positive degree. There exists a variety $Z$ such that  $\aff^1 \times Z \isom U = \bbD_+(f) \isom \Spec B_{(f)}$, 
so $B_{(f)}$ is a polynomial ring in one variable, so $f$ is a cylindrical element of $B$. 
Conversely, let $f$ be a cylindrical element of $B$ and let $U = \bbD_+(f)$.
Then (a) implies that $U$ is $H$-polar. Since $B_{(f)}$ is a polynomial ring in one variable, $U \isom \Spec B_{(f)}$ is an $H$-polar cylinder.

(c) follows immediately from (b).

(d) First note that if $n>0$ and $f \in B_n \setminus \{0\}$, then for every $D \in \Div(Y)$ 
we have $D |_{\bbD_+(f)} = 0 \Leftrightarrow D^* |_{\bbD(f)} = 0$.
Consequently, for every choice of $n>0$, $f \in B_n \setminus \{0\}$ and  $g \in (\Frac B)_0^*$, we have
\begin{multline}
\label {joiuh3ie7f8238}
H |_{ \bbD_+(f) } = \div_Y(g) |_{ \bbD_+(f) } \iff
(H-\div_Y(g)) |_{ \bbD_+(f) } = 0 \\
\iff (H-\div_Y(g))^* |_{ \bbD(f) } = 0 \iff
(\div_X(u)-\div_X(g)) |_{ \bbD(f) } = 0 \\
\iff \div_X(u/g) |_{ \bbD(f) } = 0 \iff u/g \in B_f^*,
\end{multline}
where $B_f^*$ is the set of units of $B_f$. 

Now suppose that (i) holds. Let $\pgoth \in \Proj(B)$.
Since $H$ is Cartier, there exist an open neighborhood $U$ of $\pgoth$ in $Y=\Proj(B)$ and an element $g$ of $K(Y)^* = (\Frac B)_0^*$
such that $H |_{ U } = \div_Y(g) |_{ U }$.
Choose $n>0$ and $f \in B_n \setminus \{0\}$ such that $\pgoth \in \bbD_+(f) \subseteq U$;
then $H |_{ \bbD_+(f) } = \div_Y(g) |_{ \bbD_+(f) }$.

Now suppose that (i) holds. Let $\pgoth \in \Proj(B)$.
There exist $n>0$, $f \in B_n \setminus \{0\}$ and $g \in (\Frac B)_0^*$ such that $\pgoth \in \bbD_+(f)$ and
$H |_{ \bbD_+(f) } = \div_Y(g) |_{ \bbD_+(f) }$.
By \eqref{joiuh3ie7f8238}, it follows that $u/g \in B_f^*$, so $u/g \in B_\pgoth^*$.
We claim that
\begin{equation} \label {nxbcGHFWIOJuasgdq8w2wso0k}
\text{$u/g = a/b$ for some homogeneous $a,b \in B \setminus \pgoth$.}
\end{equation}
Indeed, since $u/g \in \big( \Frac B \big)_1$,
we have $u/g = \alpha/\beta$ for some homogeneous elements $\alpha,\beta \in B \setminus \{0\}$ with $\deg(\alpha) - \deg(\beta) = 1$.
Since $u/g \in B_\pgoth^*$, we have $\alpha/\beta = \alpha'/\beta'$ for some (not necessarily homogeneous) elements $\alpha',\beta' \in B \setminus \pgoth$.
We can write $\alpha' = \sum_i \alpha'_i$ and $\beta' = \sum_i \beta'_i$ ($\alpha'_i,\beta'_i \in B_i$).
Since $\beta' \notin \pgoth$, we can choose $j$ such that $\beta'_j \notin \pgoth$.
From $\alpha \beta' = \beta \alpha'$ together with the fact that $\alpha, \beta$ are homogeneous, we deduce that
$\alpha \beta'_j = \beta \alpha'_{j+1}$, so $\alpha'_{j+1}/\beta'_j = \alpha/\beta = u/g \in B_\pgoth^*$.
If $\alpha'_{j+1} \in \pgoth$ then $\alpha'_{j+1}/\beta'_j \in \pgoth B_\pgoth$ (because $\beta_j' \notin \pgoth$), a contradiction.
So $\alpha'_{j+1} \notin \pgoth$, i.e., $a=\alpha'_{j+1}$ and $b=\beta_j'$ satisfy the requirement of \eqref{nxbcGHFWIOJuasgdq8w2wso0k}.
Since $a,b \notin \pgoth$, the degrees of $a$ and $b$ belong to $e(B/\pgoth) \Integ$;
since $\deg(a) - \deg(b) = \deg(u/g)=1$, we obtain $e(B/\pgoth)=1$.
This shows that (i) implies (ii).

Conversely, suppose that (ii) holds. Let $\pgoth \in \Proj(B)$.
Since $e(B/\pgoth)=1$, there exist homogeneous elements $a,b \in B \setminus \pgoth$ such that $\deg(a) - \deg(b) = 1$.
Then $g = ub/a$ belongs to $(\Frac B)_0^*$.
Let $f = ab$; then $\pgoth \in \bbD_+(f)$ and $u/g = a/b \in B_f^*$.
By \eqref{joiuh3ie7f8238}, it follows that $H |_{ \bbD_+(f) } = \div_Y(g) |_{ \bbD_+(f) }$.
This shows that (ii) implies (i).
\end{proof}

\begin{corollary} \label {kjxcnopo93ewjdspwqW0}
Let $\bk$ be a field, and let $B \in \GND$ and $(Y,H) \in \PP$ be such that $A(Y,H) \isom B$.
\begin{enumerata}

\item $(Y,H)$ is cylindrical if and only if $B$ has a cylindrical element.
\item $H$ is a Cartier divisor of $Y$ if and only if $e(B/\pgoth) = 1$ for all $\pgoth \in \Proj B$.

\end{enumerata}
\end{corollary}

\begin{proof}
Let $X = \Spec B$ and $Y' = \Proj B$.
Let $u \in \big( \Frac B \big)_1 \setminus \{0\}$ and let $H' \in \Div_\Rat(Y')$ be such that $(H')^* = \div_X(u)$.
Then $A(Y',H') \isom B$ by Theorem \ref{skjbci8wjqpodmpa0}, 
so $(Y,H) \approx (Y',H')$
by injectivity of the map ${\PP/\!\!\approx} \to {\GND/\!\!\isom}$\ of \ref{09203ndc020B0ej102}.

(a) By Remark \ref{xjkh9182djADBonqws7}, $(Y,H)$ is cylindrical if and only if $(Y',H')$ is cylindrical.
By Lemma \ref{oijx9w38h8f9wlodl}(c),  $(Y',H')$ is cylindrical if and only if $B$ has a cylindrical element.

(b) Since $(Y,H) \approx (Y',H')$, there exists a $\bk$-isomorphism $Y \to Y'$ that carries $H$ to a $\Rat$-divisor $D$ of $Y'$ such that $D \sim H'$.
Clearly, $H$ is a Cartier divisor of $Y$ if and only if $D$ is a Cartier divisor of $Y'$, if and only if $H'$ is Cartier.
By Lemma \ref{oijx9w38h8f9wlodl}(d), $H'$ is Cartier if and only if  $e(B/\pgoth) = 1$ for all $\pgoth \in \Proj B$.
\end{proof}

In view of part (a) of Corollary \ref{kjxcnopo93ewjdspwqW0}, it is now clear that Theorem 0.6 of \cite{KPZ2013},
as stated just before Lemma \ref{oijx9w38h8f9wlodl}, is equivalent to part (1) of assertion  \ref{iocjvn2uygew8dhc}.
We now turn our attention to:

\medskip

\noindent{\bf Corollary 3.2 of \cite{KPZ2013}. }
[{\it Let $\bk$ be an algebraically closed field of characteristic $0$.}]
{\it Let $Y$ be a normal algebraic $\bk$-variety projective over an affine variety $S$ with $\dim_S Y \ge1$.
Let $H \in \Div(Y)$ be an ample divisor on $Y$, and let $V = \Spec A(Y,H)$ be the associated affine quasicone over $Y$.
Then $V$ admits an effective $G_a$-action if and only if $Y$ contains an $H$-polar cylinder.
}

\medskip

The hypotheses of this result are equivalent to ``$(Y,H) \in \PP$ and $H$ is Cartier''.
So Corollary 3.2 of  \cite{KPZ2013} is equivalent to the following assertion:

\begin{quote}
{\it Let $\bk$ be an algebraically closed field of characteristic $0$, and let $(Y,H) \in \PP$ be such that $H$ is Cartier.
Then $A(Y,H)$ is non-rigid if and only if $(Y,H)$ is cylindrical.}
\end{quote}
So, by parts (a) and (b) of Corollary \ref{kjxcnopo93ewjdspwqW0}, 
we see that Corollary 3.2 of \cite{KPZ2013} is equivalent to part (2) of assertion \ref{iocjvn2uygew8dhc}.

\bibliographystyle{alpha}

\begin{thebibliography}{{M}at80}

\bibitem[AF88]{AndersonFeil}
M.~Anderson and T.~Feil.
\newblock {\em Lattice-ordered groups: an introduction}.
\newblock Reidel Texts in the Mathematical Sciences. D. Reidel Publishing
  Company, 1988.

\bibitem[Dai]{Dai:IntroLNDs2010}
D.~Daigle.
\newblock Introduction to locally nilpotent derivations.
\newblock Informal lecture notes prepared in 2010, available at
  \verb!http://aix1.uottawa.ca/~ddaigle!

\bibitem[Dai12]{Dai:TameWild}
D.~Daigle.
\newblock Tame and wild degree functions.
\newblock {\em Osaka J. Math.}, 49:53--80, 2012.


\bibitem[DFMJ17]{DaiFreudMoser}
Daniel Daigle, Gene Freudenburg, and Lucy Moser-Jauslin.
\newblock Locally nilpotent derivations of rings graded by an abelian group.
\newblock In {\em Algebraic varieties and automorphism groups}, volume~75 of
  {\em Adv. Stud. Pure Math.}, pages 29--48. Math. Soc. Japan, Tokyo, 2017.


\bibitem[Dem88]{Demazure_1988}
M.~Demazure.
\newblock Anneaux gradu\'{e}s normaux.
\newblock In {\em Introduction \`a la th\'{e}orie des singularit\'{e}s, {II}},
  volume~37 of {\em Travaux en Cours}, pages 35--68. Hermann, Paris, 1988.

\bibitem[Dol]{Dolg:McKayCorr}
I.~Dolgachev.
\newblock Mckay's correspondence.
\newblock Winter 2006/07, can be found at
  \verb|http://www.math.lsa.umich.edu/~idolga/lecturenotes.html|.

\bibitem[Fre17]{Freud:Book-new}
G.~Freudenburg.
\newblock {\em Algebraic theory of locally nilpotent derivations}, volume 136
  of {\em Encyclopaedia of Mathematical Sciences}.
\newblock Springer-Verlag, Berlin, second edition, 2017.
\newblock Invariant Theory and Algebraic Transformation Groups, VII.

\bibitem[Gri01]{GrilletBook_2001}
P.~A. Grillet.
\newblock {\em Commutative semigroups}, volume~2 of {\em Advances in
  Mathematics (Dordrecht)}.
\newblock Kluwer Academic Publishers, Dordrecht, 2001.

\bibitem[Har77]{Hartshorne}
R.~Hartshorne.
\newblock {\em Algebraic {G}eometry}, volume~52 of {\em {GTM}}.
\newblock Springer-{V}erlag, 1977.

\bibitem[KPZ13]{KPZ2013}
T.~Kishimoto, Y.~Prokhorov and M.~Zaidenberg.
\newblock {$\Bbb G_{\rm a}$}-actions on affine cones.
\newblock {\em Transform. Groups}, 18:1137--1153, 2013.

\bibitem[{Lan}93]{LangAlgebraText}
{S}. {Lang}.
\newblock {\em Algebra}.
\newblock Addison-Wesley, third edition, 1993.

\bibitem[MM09]{MasudaMiyanishi:Lifting09}
K. Masuda and M. Miyanishi.
\newblock Lifting of the additive group scheme actions.
\newblock {\em Tohoku Math. J. (2)}, 61(2):267--286, 2009.

\bibitem[{M}at80]{Matsumura}
{H}. {M}atsumura.
\newblock {\em {C}ommutative {A}lgebra}.
\newblock {M}athematics {L}ecture {N}ote {S}eries. {B}enjamin/{C}ummings,
  second edition, 1980.

\bibitem[NG67]{Nouaze-Gabriel_1967}
Y.~Nouaz\'{e} and P.~Gabriel.
\newblock Id\'{e}aux premiers de l'alg\`ebre enveloppante d'une alg\`ebre de
  {L}ie nilpotente.
\newblock {\em J. Algebra}, 6:77--99, 1967.

\bibitem[Vas69]{VASC}
W.~V.\ Vasconcelos.
\newblock Derivations of {C}ommutative {N}oetherian {R}ings.
\newblock {\em Math Z.}, 112:229--233, 1969.

\bibitem[Wri81]{Wright:JacConj}
D.~Wright.
\newblock On the jacobian conjecture.
\newblock {\em Illinois J.\ of Math.}, 25:423--440, 1981.



\end{thebibliography}

\end{document}